\documentclass[12pt]{article}

\usepackage{amsmath, amsfonts, amssymb}
\usepackage{graphicx}
\usepackage{epstopdf}
\usepackage{bm}
\usepackage{color}

\newcommand{\bd}[1]{\mathbf{#1}}
\newcommand{\al}{\alpha}
\newcommand{\del}{\delta}
\newcommand{\lam}{\lambda}
\newcommand{\pa}{\partial}
\newcommand{\beq}{\begin{equation}}
\newcommand{\eeq}{\end{equation}}

\DeclareMathOperator\erf{erf}
\DeclareMathOperator\erfc{erfc}

\title{Regularized Single and Double Layer Integrals\\ in 3D Stokes Flow}
\author{Svetlana Tlupova\thanks{Department of Mathematics, Farmingdale State College, SUNY, Farmingdale, NY 11735, USA tlupovs@farmingdale.edu} \and J. Thomas Beale\thanks{Department of Mathematics, Duke University, Durham, NC 27708, USA beale@math.duke.edu}}

\date{\today}

\begin{document}

\maketitle

\begin{abstract} 
We present a numerical method for computing the single layer (Stokeslet) and double layer (stresslet) integrals in Stokes flow. The method applies to smooth, closed surfaces in three dimensions, and achieves high accuracy both on and near the surface. The singular Stokeslet and stresslet kernels are regularized and, for the nearly singular case, corrections are added to reduce the regularization error. These corrections are derived analytically for both the Stokeslet and the stresslet using local asymptotic analysis. For the case of evaluating the integrals on the surface, as needed when solving integral equations, we design high order regularizations for both kernels that do not require corrections. This approach is direct in that it does not require grid refinement or special quadrature near the singularity, and therefore does not increase the computational complexity of the overall algorithm. Numerical tests demonstrate the uniform convergence rates for several surfaces in both the singular and near singular cases, as well as the importance of corrections when two surfaces are close to each other.
\end{abstract}

{\bf Keywords:} Stokes flow; Boundary integral method; Nearly singular integrals; Regularization.

\section{Introduction}

Stokes flows are of relevance in many practical problems where the length scales are very small, the fluid is viscous, or the velocity is very small (i.e., creeping flows), all resulting in a small Reynolds number. In dimensionless form, the incompressible Stokes equations are
\beq
	-\nabla p + \Delta \bd{u} = 0,\qquad \nabla\cdot \bd{u} = 0,
\eeq
where $p$ is the pressure and $\bd{u}$ is the flow velocity. The Stokeslet and stresslet are the primary fundamental solutions for the velocity $\bd{u}$:
\begin{subequations}
\begin{align}
	S_{ij}(\bd{y,x}) &= \frac{\del_{ij}}{|\bd{y} - \bd{x}|} + \frac{(y_i - x_i)(y_j - x_j)}{|\bd{y} - \bd{x}|^3}, \\[6pt]
	T_{ijk} (\bd{y,x}) &= -\frac{6(y_i - x_i)(y_j - x_j)(y_k - x_k)}{|\bd{y} - \bd{x}|^5},
\end{align}
\end{subequations}
where $\del_{ij}$ is the Kronecker delta and $i,j,k = 1,2,3$ are Cartesian coordinates.
When used in boundary integral methods, these lead to the single and double layer representations of Stokes flow, respectively:
\begin{subequations}
\begin{align}
	\label{SingleLayer}
	u_i(\bd{y}) &= \frac{1}{8\pi}\int_{\pa\Omega} S_{ij}(\bd{y,x}) f_j(\bd{x})dS(\bd{x}), \\[6pt]
	\label{DoubleLayer}
	w_i(\bd{y}) &= \frac{1}{8\pi}\int_{\pa\Omega} T_{ijk} (\bd{y,x}) q_j(\bd{x}) n_k(\bd{x})dS(\bd{x}),
\end{align}
\end{subequations}
where $n_k$ are the components of the unit outward normal vector to the surface. The integral in \eqref{SingleLayer} is continuous across $\pa\Omega$, and the integral in \eqref{DoubleLayer} is discontinuous and has a jump of $\mp 4\pi q_i(\bd{x}_0)$ in the limit from either the interior or exterior of the domain.

The importance of boundary integral equations in Stokes flow models is well recognized. In particular, they have been used extensively in interfacial dynamic simulations \cite{pozrikidis01,wang06}, such as vesicle flows \cite{veerapaneni11}, drop dynamics \cite{zinchenko05,janssen06}, and particle motion \cite{pozrikidis05}. For such applications, the jump conditions across the interface are incorporated into the integral formulation naturally, the dimensionality of the problem is reduced, and high accuracy can be achieved on the boundary and for points away from the boundary. Evaluating the integrals accurately for points near the boundary, e.g., when two interfaces are close together, is the most difficult case and is an active area of research \cite{tornberg18,qbx13,beale16,nguyen14,ying06,bruno01,helsing08,barnett15,carvalho18}. If the viscosities inside and outside are different, the interface velocity is found from an integral equation with single and double layers \cite{pozrikidis92,pozrikidis01,rahimian15,tornberg18,wang14}.

In dealing with the evaluation of nearly singular integrals, Ying, Biros and Zorin \cite{ying06} proposed an interpolation procedure. The authors considered 3D elliptic problems, with domain boundaries given by overlapping patches parameterized using coordinate charts. {In particular they computed double layer integrals for Stokes velocity as in \eqref{DoubleLayer}.} Since the errors decay rapidly away from the surface, the values near the boundary were obtained by interpolating between the values at points on the surface and points that are sufficiently separated from the surface along the surface normals. This algorithm was adapted and optimized by Sorgentone and Tornberg \cite{tornberg18} for close interactions of viscous drops with surface tension, where a spherical harmonics expansion was used to parameterize the surface. A quadrature by expansion method was developed by Kl\"{o}ckner et al. \cite{qbx13} and Barnett \cite{qbx14} for evaluation of Laplace and Helmholtz potentials, through local expansions. This method achieves exponential accuracy but requires upsampling the density on a finer grid. Af Klinteberg and Tornberg \cite{klinteberg16} applied the quadrature by expansion method to simulate spheroidal particles in periodic Stokes flow, using precomputations and the fast Ewald summation method for faster computations. A target-specific QBX method was developed by Siegel and Tornberg \cite{siegel18}, where the same accuracy is achieved using fewer terms. Another approach was followed by Bruno and Kunyansky in \cite{bruno01} for surface scattering problems, where partitions of unity were used along with an analytical resolution of the singularity by a change to polar coordinates.

The method of {Beale and coworkers} \cite{beale01,beale04,beale-lai} follows a different approach. The integrals are first regularized to remove the singularity. This introduces a regularization error, which is reduced by adding correction terms derived analytically using local expansions. Discretizing the integrals will introduce an additional error component, and discretization corrections can be derived as well to reduce this error. This work has led to further explorations for the Laplace's equation \cite{beale04,beale16}, the Stokes equations \cite{tlupova13,nguyen14}, and the Helmholtz equation \cite{nicholas}. In \cite{tlupova13}, the Stokes single layer integral, which has the Stokeslet as the kernel, was written using the gradient of the Laplacian Green's function as the kernel, and then regularization and discretization corrections were derived for the normal and tangential parts of this integral. The method of regularized Stokeslets \cite{cortez01,cortez05} is also based on desingularizing the velocity, but the approach is slightly different. It can be viewed as replacing the entire problem by a regularized one by smoothing out the forces acting on the fluid particles, whereas we choose the regularization specifically to compute the surface integrals accurately. In \cite{nguyen14}, Nguyen and Cortez derived regularization corrections for the regularized Stokeslets with an orthogonal system of coordinates, while paying special attention to preserving the incompressibility condition. 

In this paper we extend the framework described above to evaluate both the single and double layer integrals of Stokes flow in three dimensions on or near smooth, closed surfaces. We only treat the regularization component of the error, and our numerical results show that, with proper choice of the regularization parameter, high order uniform convergence {of the integrals \eqref{SingleLayer}, \eqref{DoubleLayer} is achieved, uniformly with respect to $\bd{y}$ near the surface}. For evaluation at points on the surface, as needed in solving integral equations, we design regularizations for the Stokeslet and stresslet kernels with high order accuracy, approaching $O(h^5)$. {Here $h$ is the grid spacing chosen in three-space and used in coordinate planes for the discretization of the integrals, as explained in Section 2.}  For evaluation near the surface we use simpler regularization and derive corrections to achieve accuracy uniformly about $O(h^3)$ in practice. With this approach, there are no parameters to fine-tune except the regularization parameter $\delta$. Through experimentation with this parameter we find that $\delta/h = 3$ for the case on the surface and $\delta/h = 2$ near the surface are reliable choices.  Because the correction formulas are precomputed analytically using local asymptotic analysis, the overall computational complexity of the algorithm does not increase. Also, no special gridding or quadrature is needed near the singularity, so that the spacing does not change with proximity to the boundary. Our method here is more direct than in \cite{tlupova13} and treats the double layer or stresslet integral as well as the Stokeslet. As in \cite{beale16}, we use a quadrature rule for closed surfaces, introduced in \cite{wilson}, which works well for general surfaces without requiring coordinate charts.

The numerical method is summarized in Section 2, including the formulas for the regularization corrections used for points near the surface.  These formulas are derived analytically in Sections 3.1 and 3.2.  The higher order regularizations for evaluation on the surface are given in Section 3.3. Various numerical examples are presented in Section 4 which illustrate the predicted performance.  Examples with known exact solutions are used to test the single layer and double layer separately (Sect. 4.2 and 4.3) and in combination (Sect. 4.4).  The integral equation for the velocity of an interface with surface tension separating two different viscosities is solved for several surfaces (Sect. 4.5) and for two spheres that are close to each other (Sect. 4.6).


\section{Numerical Method}

At the heart of the numerical method presented here is the regularization of the singularities that develop in the kernels of \eqref{SingleLayer}-\eqref{DoubleLayer} as $r=|\bd{y} - \bd{x}|$ approaches $0$. As a first step, we use subtraction to reduce the singularity in the double layer \eqref{DoubleLayer},
\begin{equation}
	\label{DoubleLayer_subtract}
	w_i(\bd{y}) = \frac{1}{8\pi}\int_{\pa\Omega} T_{ijk} (\bd{y,x}) [q_j(\bd{x}) - q_j(\bd{x}_0)] n_k(\bd{x}) dS(\bd{x}) + \frac{1}{8\pi}\chi (\bd{y}) q_i(\bd{x}_0),
\end{equation}
where $\bd{x}_0$ is the boundary point closest to $\bd{y}$, and we have applied the well known identity (see, for example, \cite{pozrikidis92} sec. 2.1-2.3, or \cite{pozrikidis01})
\begin{equation}
	\int_{\pa\Omega} T_{ijk} (\bd{y,x}) n_k(\bd{x})dS(\bd{x}) = \chi (\bd{y})\del_{ij},
\end{equation}
where $\chi (\bd{y}) = 8\pi, 4\pi, 0$ if $\bd{y}$ is inside, on, and outside the boundary, respectively. 

We then regularize the Stokeslet as follows
\begin{equation}
	\label{StokesletHH}
	\bd{u}_\del(\bd{y}) = \frac{1}{8\pi}\int_{\pa\Omega} \left[\bd{f}(\bd{x}) \frac{s_1(r/\del)}{r} + (\bd{f}(\bd{x})\cdot(\bd{y}-\bd{x}))(\bd{y}-\bd{x}) \frac{s_2(r/\del)}{r^3} \right] dS(\bd{x}), 
\end{equation}
where the smoothing factors are chosen so that $\lim_{\rho\to\infty}s_1(\rho)=1$, $\lim_{\rho\to\infty}s_2(\rho)=1$, $s_1(\rho)=O(\rho)$ and $s_2(\rho)=O(\rho^3)$ for $\rho$ small, and $s_1(r/\del)/r$, $s_2(r/\del)/r^3$ are smooth as functions of $\bf y - \bf x$ with $r = |\bf y - \bf x|$ and fixed $\delta>0$.

The regularized version of the stresslet is similar,
\begin{align}
	\label{StressletHH}
	\bd{w}_\del(\bd{y}) = -\frac{3}{4\pi} \int_{\pa\Omega} &[(\bd{y}-\bd{x})\cdot \tilde{\bd{q}}(\bd{x})] [(\bd{y}-\bd{x})\cdot \bd{n}(\bd{x})] (\bd{y}-\bd{x}) \frac{s_3(r/\del)}{r^5} dS(\bd{x})  + \frac{1}{8\pi}\chi (\bd{y}) \bd{q}(\bd{x}_0), 
\end{align}
where $\tilde{\bd{q}}(\bd{x}) = \bd{q}(\bd{x})-\bd{q}(\bd{x}_0)$, and 
$s_3$ is chosen with $\lim_{\rho\to\infty}s_3(\rho)=1$, $s_3(\rho)=O(\rho^5)$ for small $\rho$, and $s_3(r/\del)/r^5$ smooth for $\delta>0$.

There are a number of possible choices for the smoothing factors; several were explored in \cite{nguyen14}. In this work, we will use
\begin{subequations}
\begin{align}
	\label{s1}
	s_1(r) &= \erf(r), \\
	s_2(r) &= \erf(r) - 2 r e^{-r^2}/\sqrt{\pi}, \label{s2}\\
	\label{s3}
	s_3(r) &= \erf(r) - 2 r \left(\frac23 r^2 + 1\right) e^{-r^2}/\sqrt{\pi},
\end{align}
\end{subequations}
where $\erf$ is the error function. This choice of regularization is simple; $(s_1-1)/r$ decays rapidly in the far field, and $s_2, s_3$ are derived by modifying the error function to get $s_2(\rho)=O(\rho^3)$ and $s_3(\rho)=O(\rho^5)$ for small $\rho$ as mentioned above.

The regularization error, defined as $\bd{u}_\del - \bd{u}$ for the Stokeslet and $\bd{w}_\del - \bd{w}$ for the stresslet, is at best $O(\del)$. Our approach to reducing this error is twofold. When the integrals are evaluated on the surface (e.g., when solving integral equations), we modify the smoothing functions \eqref{s1}-\eqref{s3} to readily achieve a higher accuracy of $O(\del^5)$. This is explained in Section 3.3. When the evaluation point is near the surface, however, such a modification is not helpful, and we increase the accuracy by adding corrections. These corrections represent the dominant terms in the regularization error, and improve the accuracy to $O(\del^3)$. Specifically, we compute the single layer and the double layer as 
\begin{subequations}
\begin{align}
	\label{SL_tilde}
	\tilde{\bd{u}}_\del(\bd{y}) &= \bd{u}_\del(\bd{y}) + \bd{C}^\bd{u}(\bd{x}_0), \\
	\label{DL_tilde}
	\tilde{\bd{w}}_\del(\bd{y}) &= \bd{w}_\del(\bd{y}) + \bd{C}^\bd{w}(\bd{x}_0), 
\end{align}
\end{subequations}
with corrections, to be explained,
\begin{subequations}
\begin{align}
	\label{Corr_SL}
    	\bd{C}^\bd{u}(\bd{x}_0) = &- \frac{\del}{8} \Big\{(1+ H\lam\del)
             [ 2 (I_1+I_{2a}) \bd{f}^{\rm nl} + (2I_1+I_{2b}) \bd{f}^{\rm tan} ]
         - \del \lam I_{2b} [\nabla_S(\bd{f}\cdot \bd{n}) + (\nabla_S\cdot\bd{f})\bd{n}] \Big\},\\
	\bd{C}^\bd{w}(\bd{x}_0) = &-\frac34 \del (1+\lambda\del H)I_{3a}
         \Big[ \nabla_S(\tilde{\bd{q}}\cdot\bd{n}) + (\nabla_S\cdot \tilde{\bd{q}})\bd{n} \Big]
         + \frac38 \del^2\lambda I_{3a} [\Delta_S \tilde{\bd{q}}^{\rm nl}]   \nonumber \\
         \label{Corr_DL}
         &+ \frac{3}{32}\del^2\lambda I_{3b} \Big[
               \{ \Delta_S (\tilde{\bd{q}} + \tilde{\bd{q}}^{\rm nl}) \}^{\rm tan}
           + 2\nabla_S( \nabla_S\cdot \tilde{\bd{q}} ) - 4H \nabla_S(\tilde{\bd{q}}\cdot\bd{n})\Big].
\end{align}
\end{subequations}
Here, assuming that $\bd{y}$ is near the surface, $\bd{x}_0$ is the surface point closest to $\bd{y}$, so that $\bd{y}=\bd{x}_0+b\bd{n}$ for some $b$, where $\bd{n}$ is the unit outward normal at $\bd{x}_0$.
$H$ is the mean curvature and $\lam = b/\del$.  The superscript ${{\rm nl}}$ denotes the normal part of a vector, $\bd{a}^{{\rm nl}} \equiv (\bd{a}\cdot \bd{n}) \bd{n}$, and superscript ${{\rm tan}}$ denotes the tangential part,
$\bd{a}^{{\rm tan}} \equiv \bd{a}-\bd{a}^{{\rm nl}}$.
For any local coordinate system $\bd{x}=\bd{x}(\al_1,\al_2)$ on the surface, we have
tangent vectors $\bd{T}_i = \pa\bd{x}/\pa \al_i$, dual vectors $\bd{T}_i^*$,
$i = 1,2$, such that $\bd{T}_i^*\cdot\bd{T}_j = \delta_{ij}$, and the metric tensor
$g_{ij} =  \bd{T}_i\cdot\bd{T}_j$.  Note that $\bd{a}^{{\rm tan}} = (\bd{a}\cdot\bd{T}_1)\bd{T}_1^* + 
(\bd{a}\cdot\bd{T}_2)\bd{T}_2^*$.
The surface gradient of a scalar function $f$, the surface divergence of a vector function $\bd{v}$, and the surface Laplacian of a scalar function $f$ are defined as
\begin{subequations}
\begin{align}
	\label{surf_gradient}
	\nabla_S f &= \frac{\pa f}{\pa \al_1} \bd{T}_1^* + \frac{\pa f}{\pa \al_2} \bd{T}_2^*, \\
	\label{surf_divergence}
	\nabla_S\cdot\bd{v} &= \frac{\pa \bd{v}^{{\rm tan}}}{\pa \al_1} \cdot \bd{T}_1^* + \frac{\pa \bd{v}^{{\rm tan}}}{\pa \al_2} \cdot \bd{T}_2^*, \\
	\label{surf_laplacian}
	\Delta_S f &= \frac{1}{\sqrt{g}}\sum_{i,j=1}^2 \frac{\pa}{\pa \al_j}\Big( \sqrt{g}g^{ij}\frac{\pa f}{\pa\al_i}\Big).
\end{align}
\end{subequations}
They are independent of coordinates.  It can be checked that (\ref{surf_divergence})
agrees with formula (9.41.1) in \cite{aris} for the surface divergence of a tangential vector field.
Section 4.1 has some details on how these quantities are computed for a Monge parameterization. {In \eqref{Corr_SL}, \eqref{Corr_DL} $I_1,\dots,I_{3b}$ are certain integrals that occur in the derivations of Section 3.}  With our choice of $s_1, s_2, s_3$,
\begin{subequations}
\begin{align}
	I_1(\lam) &= |\lam| \erfc|\lam|  - e^{\lam^2}/\sqrt{\pi} \\
	I_{2a}(\lam) &= -|\lam| \erfc|\lam| \\
	I_{2b}(\lam) &= 2 I_1 \\
	I_{3a}(\lam) &= -\frac23 |\lam| \erfc|\lam| \\
	I_{3b}(\lam) &= \frac83 \left( |\lam| \erfc|\lam|  - e^{\lam^2}/\sqrt{\pi} \right)
\end{align}
\end{subequations}
where $\erfc(x)=1-\erf(x)$.  If we chose (\ref{s1},\ref{s2},\ref{s3}) differently, the expressions
for the corrections would be the same, but $I_1,\dots,I_{3b}$ would be different.

Once the integrands are smoothed out, we discretize the integrals using a surprisingly simple quadrature method for closed surfaces introduced in \cite{wilson} and explained in \cite{beale16}.  We choose an angle $\theta$ and define a partition of unity on the unit sphere,
consisting of functions $\psi_1, \psi_2, \psi_3$ with $\Sigma_i \psi_i \equiv 1$ such that
$\psi_i(\bd{n}) = 0$ if $|\bd{n}\cdot\bd{e}_i| \leq \cos{\theta}$, where $\bd{e}_i$ is the
$i$th coordinate vector.   Here we choose $\theta = 70^o$.  For mesh size $h$,
a set $R_3$ of quadrature points consists of points $\bf{x}$ on the surface of the form
$(j_1h,j_2h,x_3)$ such that $|\bd{n(x}) \cdot \bd{e}_3| \geq \cos{\theta}$, where
$\bd{n}(\bd{x})$ is the unit normal at $\bd{x}${, see Fig.~\ref{Quad_pts}}.  Sets $R_1$ and $R_2$ are
defined similarly.  For a function $f$ on the surface the integral is computed as
\beq
	\label{Quadrature}
	\int_S f(\bd{x}) \,dS(\bd{x}) \approx \sum_{i=1}^3 \sum_{\bd{x} \in R_i} 
            \frac{\psi_i(\bd{n}(\bd{x}))\,f(\bd{x}) }{| \bd{n}(\bd{x})\cdot\bd{e}_i |} \,h^2
\eeq
{The partition of unity functions $\psi_i$ are constructed from the $C^\infty$ bump function $b(r) = \exp{((r^2/(r^2 - 1))}$ for $|r| < 1$ and zero otherwise.  The quadrature is effectively reduced to the trapezoidal rule without boundary.  Thus for regular integrands the quadrature has arbitrarily high order accuracy, limited only by the degree of smoothness of the integrand and surface.}  The points in $R_i$
can be found by a line search since they are well separated; see \cite{wilson} and \cite{beale16}.
\begin{figure}[htb]
\centering
\includegraphics[scale=0.5]{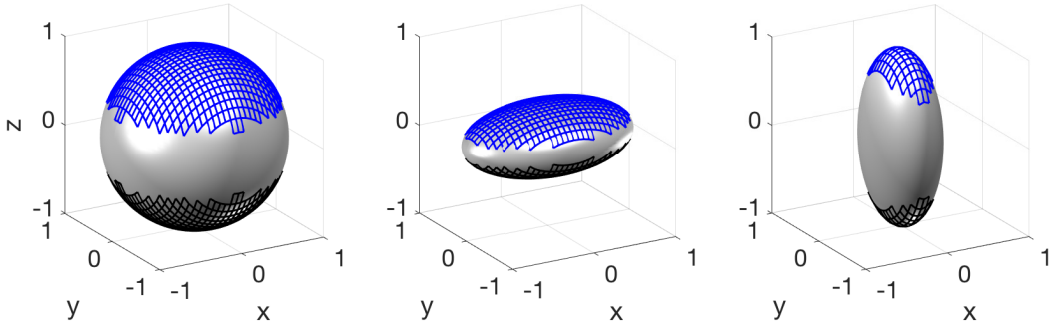}
\caption{Quadrature points generated by $\psi_3$ and $h=1/16$ on the surface of the unit sphere (left), the ellipsoid with semiaxes $a=1,b=0.6,c=0.4$ (middle), and the ellipsoid $a=0.4,b=0.6,c=1$ (right). The quadrature points are at the intersections of the lines.}
\label{Quad_pts}
\end{figure}

Finally we discuss the accuracy of this method, relying on the analogy with the case of Laplace's equation treated in \cite{beale04}, \cite{beale16}.  Error estimates for the harmonic double layer potential were derived in \cite{beale04} and extended to the single layer in \cite{beale16}.  The present work is similar except that discretization corrections were included in the earlier case.  Based on the theory for the harmonic potentials, we expect that the error in the present method, for evaluation of (\ref{SingleLayer}, \ref{DoubleLayer}) at points near the surface, can be estimated uniformly as
$$  \epsilon \leq C_1\delta^3 + C_2 h\,e^{-c_0(\delta/h)^2}  $$
The first term represents the regularization error remaining after the correction
(\ref{Corr_SL}, \ref{Corr_DL}), and the second term is the discretization error.  For evaluation on the surface, using the formulas of Sec. 3.3, the error estimates for the single and double layer integrals should be
$$  \epsilon_u \leq C_1\delta^5 + C_2 h\,e^{-c_0(\delta/h)^2}  \, \qquad
   \epsilon_w \leq C_1\delta^5 + C_2 h^2\,e^{-c_0(\delta/h)^2}  $$
The improvement in the double layer results from the subtraction {in equation \eqref{DoubleLayer_subtract}}.  These estimates are discussed further in Sec. 3.4 of \cite{beale16}.  The accuracy depends critically on the relationship between $\delta$ and $h$.  In this work we have taken
$\delta/h =$ constant, such as $2$ or $3$, for simplicity.  {In the numerical examples reported in Section 4, we see errors about $O(h^3)$ near the surface and significantly higher order on the surface, indicating that the regularization error is dominant.  However, this pattern could not continue as $h \to 0$ with $\delta/h$ fixed.  In principle the discretization error can be controlled as $h \to 0$ by increasing $\delta/h$, e.g. $\delta = ch^p$ for any chosen $p < 1$.  Then the exponential in the error estimates decreases rapidly, and the second term, representing the discretization error, is dominated by the first term as $h \to 0$.  In this way convergence can be achieved which is $O(h^{3p})$ near the surface and higher on the surface.}


\section{Regularization Corrections}

In this section, we derive the regularization corrections for the Stokeslet and stresslet integrals, and discuss the special case of evaluating them on the surface with very high accuracy.

\subsection{Stokeslet}

Computing the regularization correction in the first part of the Stokeslet \eqref{StokesletHH} is similar to the derivation for the Laplace's single layer potential in \cite{beale16}, but the second part needs extra care. First write the error as
\begin{equation}
	\epsilon = \frac{1}{8\pi}\int_{\pa\Omega} \Big\{ \bd{f}(\bd{x}) \frac{\phi_1(r/\del)}{r} + (\bd{f}(\bd{x})\cdot(\bd{y}-\bd{x}))(\bd{y}-\bd{x}) \frac{\phi_2(r/\del)}{r^3} \Big\}dS(\bd{x}),
\end{equation}
where we define $\phi(r)=s(r)-1$.  For the part of the surface near $\bd{y}$ we will use a special parameterization $\bd{x}(\alpha)$, with $\bd{x}(0)=\bd{0}$ and 
$\bd{y}$ along the normal line from $\bd{x}(0)$, so that $\bd{y}=b\bd{n}_0$ for some $b$, where $\bd{n}_0$ is the unit normal at $\bd{x}(0)$.  With tangent vectors to the surface $\bd{T}_i = \pa\bd{x}/\pa \al_i$, $i=1,2$ and metric tensor $g_{ij} = \bd{T}_i\cdot\bd{T}_j$, we assume that, at $\alpha = 0$, $g_{ij} = \delta_{ij}$ and
$\pa g_{ij}/\pa\al_k = 0$, $i,j,k=1,2$. Also, rotating if necessary, $\bd{T}_1, \bd{T}_2$ have the directions of principal curvature, and $\bd{n}_0=\bd{T}_1\times\bd{T}_2$.  We use expansions near $\al = 0$,
\begin{subequations}
\begin{align}
        \label{expandx}
	\bd{x}(\al) &= \bd{T}_i\al_i + \frac{1}{2}\kappa_i \bd{n}_0 \al_i^2 + O(|\al|^3),\\
	\label{expandn}
        \bd{n}(\al) &= \bd{n}_0 - \kappa_i\bd{T}_i\al_i + O(|\al|^2), \\
	\bd{f}(\al) &= \bd{f}_0+\bd{f}_i \al_i + O(|\al|^2),
\end{align}
\end{subequations}
where $\bd{T}_i = \bd{T}_i(0)$,
$\bd{f}_0 = \bd{f}(0)$, $\bd{f}_i=\pa\bd{f}/\pa\al_i (0)$, and summation in $i=1,2$ is assumed. Then 
\begin{align}
    	\bd{f}(\al) \cdot (\bd{y}-\bd{x}(\al)) &= (\bd{f}_0\cdot\bd{n}_0)b + (\bd{f}_i\cdot\bd{n}_0)b\al_i  - (\bd{f}_0\cdot\bd{T}_i)\al_i \nonumber \\
    	&- \frac{1}{2}\kappa_i(\bd{f}_0\cdot\bd{n}_0)\al_i^2 - (\bd{f}_i\cdot\bd{T}_j)\al_i\al_j + O(|\al|^3).
\end{align}
In the expansion of $\bd{F}:=[\bd{f}(\al) \cdot (\bd{y}-\bd{x}(\al))](\bd{y}-\bd{x}(\al))$, only the terms even in $\al$ will contribute:
\begin{align}
    	\bd{F}^{even} &= (\bd{f}_0\cdot\bd{n}_0)\bd{n}_0 b^2 - (\bd{f}_0\cdot\bd{n}_0)\bd{n}_0 \kappa_i b\al_i^2 + (\bd{f}_0\cdot\bd{T}_i)\bd{T}_j \al_i\al_j \nonumber\\
    	& - (\bd{f}_i\cdot\bd{n}_0)\bd{T}_j b\al_i\al_j - (\bd{f}_i\cdot\bd{T}_j)\bd{n}_0 b\al_i\al_j + O(|\al|^4 + |\al|^3b).
\end{align}
We choose a new parameter $\xi$ to replace $\al$, defined by $r^2 = b^2 + |\xi|^2$ and $\xi_i/|\xi| = \al_i/|\al|$, 
\begin{equation}
    	\al_i = (1+\frac{b\mu}{2})\xi_i + O(|\xi|^3 + b^3),
\end{equation}
\begin{equation}
    	\mu = \kappa_1\frac{\xi_1^2}{|\xi|^2} + \kappa_2\frac{\xi_2^2}{|\xi|^2}.
\end{equation}
\begin{align}
    	\bd{F}^{even}_1 &= (\bd{f}_0\cdot\bd{n}_0)\bd{n}_0 b^2 - (\bd{f}_0\cdot\bd{n}_0)\bd{n}_0 \kappa_i b\xi_i^2 + (\bd{f}_0\cdot\bd{T}_i)\bd{T}_i (1+b\mu)\xi_i^2 \nonumber\\
    	& - (\bd{f}_i\cdot\bd{n}_0)\bd{T}_i b\xi_i^2 - (\bd{f}_i\cdot\bd{T}_i)\bd{n}_0 b\xi_i^2 + O(|\xi|^4 + b^4).
\end{align}
Note that $\bd{F}_1^{even}$ is $\bd{F}^{even}$ where terms with $\xi_1\xi_2$ were omitted, as they lead to $\int_0^{2\pi} \cos(\theta)\sin(\theta)d\theta$ or $\int_0^{2\pi} \cos^3(\theta)\sin(\theta)d\theta$ type coefficients once we change $\xi$ to polar coordinates, and thus contribute zero.
The regularization error is now written as
\begin{equation}
	\label{epsilon_Reg}
    	\epsilon = \frac{1}{8\pi}\int \Big\{ \bd{m}_1(\xi,b) \frac{\phi_1\big(\sqrt{|\xi|^2+b^2}/\del\big)}{(|\xi|^2+b^2)^{1/2}} + \bd{m}_2(\xi,b) \frac{\phi_2\big(\sqrt{|\xi|^2+b^2}/\del\big)}{(|\xi|^2+b^2)^{3/2}} \Big\} d\xi,
\end{equation}
where $\bd{m}_1$ includes the nonradial terms,
\begin{equation}
    	\bd{m}_1(\xi,b) = \bd{f}\Big|\frac{\pa\al}{\pa\xi}\Big| |\bd{T}_1\times\bd{T}_2|,
\end{equation}
where $|\pa\al/\pa\xi| = 1+b\mu + O(|\xi|^2+b^2)$, and $|\bd{T}_1\times\bd{T}_2| = 1+O(|\xi|^2)$. Combining the terms in $\bd{m}_1$ and neglecting the terms odd in $\xi$ (since they contribute 0 to the integral in $\xi$), we get 
\begin{equation}
	\label{w1_even}
    	\bd{m}_1^{even}(\xi,b) = \bd{f}_0(1+b\mu) + O(|\xi|^2 + b^2).
\end{equation}
Similarly for $\bd{m}_2$,
\begin{align}
    	\bd{m}_2^{even}(\xi,b) &= \bd{F}^{even}_1\Big|\frac{\pa\al}{\pa\xi}\Big| |\bd{T}_1\times\bd{T}_2| \nonumber\\
    	&= (\bd{f}_0\cdot\bd{n}_0)\bd{n}_0 b^2 (1+b\mu) - (\bd{f}_0\cdot\bd{n}_0)\bd{n}_0 b\mu |\xi|^2 + (\bd{f}_0\cdot\bd{T}_i)\bd{T}_i (1+2b\mu)\xi_i^2 \nonumber\\
	\label{w2_even}
    	& - (\bd{f}_i\cdot\bd{n}_0)\bd{T}_i b\xi_i^2 - (\bd{f}_i\cdot\bd{T}_i)\bd{n}_0 b\xi_i^2 + O(|\xi|^4 + b^4).
\end{align}

Substituting \eqref{w1_even} and \eqref{w2_even} into \eqref{epsilon_Reg}, and changing to polar coordinates such that $\xi = \del\zeta$, $b=\del\lambda$, $|\zeta| = \eta$, $q = \kappa_1 \cos^2\theta + \kappa_2 \sin^2\theta$, we get
\begin{equation}
	\label{eps_eta}
    	\epsilon = \frac{1}{8\pi} \int_0^{2\pi}\int_0^{\infty} \Big\{ \bd{m}_1^{even}(\del\zeta,\del\lambda)  \frac{\del \phi_1(\sqrt{\eta^2+\lambda^2})}{(\eta^2+\lambda^2)^{1/2}} + \bd{m}_2^{even}(\del\zeta,\del\lambda) \frac{\phi_2(\sqrt{\eta^2+\lambda^2})}{\del (\eta^2+\lambda^2)^{3/2}} \Big\} \eta d\eta d\theta,
\end{equation}
which simplifies to
\begin{subequations}
\begin{align}
	\label{e_1}
    	\epsilon = \frac{\del}{8} & \left\{2\bd{f}_0 (1+H\lambda\del) I_1 + 2(\bd{f}_0\cdot\bd{n}_0)\bd{n}_0(1+ H\lambda\del)I_{2a} - 2(\bd{f}_0\cdot\bd{n}_0)\bd{n}_0 H \lambda\del I_{2b}\right. \\
	\label{e_2}
    	&+ (\bd{f}_0\cdot\bd{T}_i)\bd{T}_i I_{2b} + (\bd{f}_0\cdot\bd{T}_i)\bd{T}_i (H+\kappa_i)\lambda\del I_{2b} \\
	\label{e_last}
    	&- \left.\left[ (\bd{f}_i\cdot\bd{n}_0)\bd{T}_i + (\bd{f}_i\cdot\bd{T}_i)\bd{n}_0 \right] \lambda\del I_{2b} \right\} + O(\del^3),
\end{align}
\end{subequations}
where $H = (\kappa_1 + \kappa_2)/2$ is the mean curvature, and
\begin{subequations}
\begin{align}
	\label{I1}
   	I_1 =& \int_0^{\infty} \frac{\phi_1(\sqrt{\eta^2+\lambda^2})}{(\eta^2+\lambda^2)^{1/2}} \eta d\eta, \\
	\label{I2}
   	I_{2a} =& \int_0^{\infty} \frac{\phi_2(\sqrt{\eta^2+\lambda^2})}{(\eta^2+\lambda^2)^{3/2}} \lambda^2\eta d\eta, \\
    	\label{I3}
    	I_{2b} =& \int_0^{\infty} \frac{\phi_2(\sqrt{\eta^2+\lambda^2})}{(\eta^2+\lambda^2)^{3/2}} \eta^3 d\eta.
\end{align}
\end{subequations}

This gives the correction expressed in a special coordinate system. To extend this to an arbitrary system, first we find from \eqref{expandx}-\eqref{expandn} that in our special coordinates at $\al = 0$, we have 
$\pa_i \bd{T}_i = \kappa_i \bd{n}_0$ and $\pa_i \bd{n} = - \kappa_i \bd{T}_i$, $i = 1,2$, where $\pa_i$ is the partial derivative in $\al_i$.  Using these we get at $\al = 0$
\begin{subequations}
\begin{align}
	\label{SurfGrad}
	&[\pa_i(\bd{f} \cdot \bd{n}) ] \bd{T}_i = (\bd{f}_i \cdot \bd{n}_0) \bd{T}_i  - \kappa_i(\bd{f} \cdot \bd{T}_i)\bd{T}_i, \\
	\label{SurfDiv}
	&[\pa_i(\bd{f} \cdot \bd{T}_i)] \bd{n} = (\bd{f}_i \cdot \bd{T}_i) \bd{n}_0  + \kappa_i(\bd{f} \cdot \bd{n}_0)\bd{n}_0.
\end{align}
\end{subequations}
Multiplying by $-1$, we see that the two terms on the right side of \eqref{SurfGrad} occur in
\eqref{e_2} and \eqref{e_last}. Similarly the two terms in \eqref{SurfDiv} occur in \eqref{e_1} and \eqref{e_last}.  We can combine terms to get, with the sum over $i=1,2$ implied,
\begin{align}
    	\epsilon = \frac{\del}{8} & \Big\{2\bd{f}_0 (1+H\lambda\del) I_1 + 2(\bd{f}_0\cdot\bd{n}_0)\bd{n}_0(1+ H\lambda\del)I_{2a} + (\bd{f}_0\cdot\bd{T}_i)\bd{T}_i (1+ H\lambda\del) I_{2b} \nonumber\\
	\label{e_final}
    	& - [\pa_i(\bd{f} \cdot \bd{n})]_0 \bd{T}_i \lambda\del I_{2b} - [\pa_i(\bd{f} \cdot \bd{T}_i)]_0 \bd{n}_0 \lambda\del I_{2b} \Big\} + O(\del^3).
\end{align}
The dual basis in the tangent space is defined as
$\bd{T}_i^* = \Sigma_j g^{ij}\bd{T}_j$,
so that $\bd{T}_i^*\cdot \bd{T}_j = \del_{ij}$.
Thus $\bd{T}_i^* = \bd{T}_i$ at $\al = 0$,
and $\Sigma_i (\bd{f}_0\cdot\bd{T}_i)\bd{T}_i =\bd{f}_0^{{\rm tan}}$.  Similarly
we can replace $\bd{T}_i$ with
$\bd{T}_i^*$ in  $[\pa_i(\bd{f} \cdot \bd{n})]_0 \bd{T}_i$, and the sum over $i$ is
$\nabla_S(\bd{f}\cdot \bd{n})$, the surface gradient of the scalar function
$\bd{f}\cdot \bd{n}$.  Also $\Sigma_i \pa_i (\bd{f} \cdot \bd{T}_i) =
\Sigma_i \pa_i (\bd{f}^{{\rm tan}} \cdot \bd{T}_i) = 
\Sigma_i \pa_i \bd{f}^{{\rm tan}} \cdot \bd{T}_i =
\Sigma_i \pa_i \bd{f}^{{\rm tan}} \cdot \bd{T}_i^* =
\nabla_S\cdot\bd{f} $,
the surface divergence at $\al = 0$, 
according to \eqref{surf_divergence}; we used the fact that
$\pa_i\bd{T}_i$ is normal.
We can now write the regularization error as
\begin{equation}
\epsilon = \frac{\del}{8} \Big\{(1+ H\lambda\del)
             [ 2 I_1 \bd{f}_0 + 2 I_{2a} \bd{f}_0^{{\rm nl}} + I_{2b} \bd{f}_0^{{\rm tan}} ]
         - \del \lambda I_{2b} [\nabla_S(\bd{f}\cdot \bd{n}) + (\nabla_S\cdot\bd{f})\bd{n}_0]
                                 \Big\}  + O(\del^3),
\end{equation}
where $\bd{f}_0^{{\rm nl}}$ and $\bd{f}_0^{{\rm tan}}$ are the normal and tangential parts of $\bd{f}_0$, leading to the correction \eqref{Corr_SL}.


\subsection{Stresslet}

The regularization error in the stresslet can be written as
\begin{equation}
	\label{e_stress}
	\epsilon = -\frac{3}{4\pi}\int_{\pa\Omega} [(\bd{y}-\bd{x})\cdot\tilde{\bd{q}}(\bd{x})] [(\bd{y}-\bd{x})\cdot \bd{n}] (\bd{y}-\bd{x}) \frac{\phi_3(r/\del)}{r^5} dS(\bd{x}),
\end{equation}
where as before, $\tilde{\bd{q}}(\bd{x}) = \bd{q}(\bd{x})-\bd{q}(\bd{x}_0)$. In what follows, we will write $\bd{q}$ for $\tilde{\bd{q}}$ and assume $\bd{q}(0)=0$ because of the subtraction. We compute the Taylor expansions of the nonradial parts:
\begin{subequations}
\begin{align}
	\label{y_minus_x_dot_n}
    	(\bd{y}-\bd{x}(\al))\cdot \bd{n}(\al) &= b + \frac{1}{2}\kappa_i\al_i^2 + O(|\al|^3 + b|\al|^2),\\[6pt]
	\label{y_minus_x_dot_g}
    	(\bd{y}-\bd{x}(\al))\cdot \bd{q}(\al) &= (\bd{q}_i\cdot\bd{n}_0)b\al_i + \frac{1}{2}(\bd{q}_{ij}\cdot\bd{n}_0)b\al_i\al_j - \frac{1}{2}(\bd{q}_i\cdot\bd{n}_0)\al_i\kappa_j\al_j^2 \nonumber \\
    	& - (\bd{q}_i\cdot\bd{T}_j)\al_i\al_j - \frac{1}{2}(\bd{q}_{ij}\cdot\bd{T}_l)\al_i\al_j\al_l + O(|\al|^4 + b|\al|^3).
\end{align}
\end{subequations}
Next we compute the product of \eqref{y_minus_x_dot_n} and \eqref{y_minus_x_dot_g}:
\begin{align}
    	[(\bd{y}-&\bd{x}(\al))\cdot \bd{q}(\al)] [(\bd{y}-\bd{x}(\al))\cdot \bd{n}(\al)] = (\bd{q}_i\cdot\bd{n}_0)b^2\al_i + \frac{1}{2}(\bd{q}_{ij}\cdot\bd{n}_0)b^2\al_i\al_j \nonumber \\
    	& - (\bd{q}_i\cdot\bd{T}_j)b\al_i\al_j - \frac{1}{2}(\bd{q}_{ij}\cdot\bd{T}_l)b\al_i\al_j\al_l - \frac{1}{2}(\bd{q}_i\cdot\bd{T}_j)\kappa_l\al_i\al_j\al_l^2 + O(|\al|^5 + b^5).
\end{align}
Now, let 
\begin{align}
    	\bd{Q} = [(\bd{y}-\bd{x}(\al))\cdot \bd{q}(\al)] [(\bd{y}-\bd{x}(\al))\cdot \bd{n}(\al)] \ (\bd{y}-\bd{x}(\al)).
\end{align}
In the expansion of this quantity, we keep only terms even in $\al$, and we also neglect terms with $\al_1\al_2$, $\al_1^3\al_2$, and $\al_1\al_2^3$, as all of these will contribute 0 to the error. We then get the following:
\begin{align}
	\bd{Q}^{even} = &- (\bd{q}_i\cdot\bd{T}_i)\bd{n}_0 b^2\al_i^2 - (\bd{q}_i\cdot\bd{n}_0)\bd{T}_ib^2\al_i^2 \nonumber\\
	&+ \frac{1}{2}(\bd{q}_{ii}\cdot\bd{n}_0)\bd{n}_0b^3\al_i^2 + \frac{1}{2}(\bd{q}_{11}\cdot\bd{T}_1)\bd{T}_1b\al_1^4 + \frac{1}{2}(\bd{q}_{22}\cdot\bd{T}_2)\bd{T}_2b\al_2^4 \nonumber\\
	&+ \Big[ \frac{1}{2}(\bd{q}_{22}\cdot\bd{T}_1)\bd{T}_1 + (\bd{q}_{12}\cdot\bd{T}_2)\bd{T}_1 + (\bd{q}_{12}\cdot\bd{T}_1)\bd{T}_2 + \frac{1}{2}(\bd{q}_{11}\cdot\bd{T}_2)\bd{T}_2 \Big] b\al_1^2 \al_2^2 \nonumber\\
	&+ O(|\al|^6 + b^6).
\end{align}
Next, we rewrite $\bd{Q}^{even}$ in the new parameter $\xi$, and then write the regularization error \eqref{e_stress} as
\begin{equation}
	\label{e_stress_1}
    	\epsilon = -\frac{3}{4\pi} \int \bd{m}(\xi,b) \frac{\phi_3\big(\sqrt{|\xi|^2+b^2}/\del\big)}{(|\xi|^2+b^2)^{5/2}} d\xi,
\end{equation}
where
\begin{align}
	\label{w_stress}
    	\bd{m}(\xi,b) = & \bd{Q}^{even} \Big|\frac{\pa\al}{\pa\xi}\Big| |\bd{T}_1\times\bd{T}_2| \nonumber\\
	= &- (\bd{q}_i\cdot\bd{T}_i)\bd{n}_0 b^2 (1+2b\mu)\xi_i^2 - (\bd{q}_i\cdot\bd{n}_0)\bd{T}_ib^2 (1+2b\mu)\xi_i^2 \nonumber\\
	&+ \frac{1}{2}(\bd{q}_{ii}\cdot\bd{n}_0)\bd{n}_0b^3\xi_i^2 + \frac{1}{2}(\bd{q}_{11}\cdot\bd{T}_1)\bd{T}_1b\xi_1^4 + \frac{1}{2}(\bd{q}_{22}\cdot\bd{T}_2)\bd{T}_2b\xi_2^4 \nonumber\\
	&+ \Big[ \frac{1}{2}(\bd{q}_{22}\cdot\bd{T}_1)\bd{T}_1 + (\bd{q}_{12}\cdot\bd{T}_2)\bd{T}_1 + (\bd{q}_{12}\cdot\bd{T}_1)\bd{T}_2 + \frac{1}{2}(\bd{q}_{11}\cdot\bd{T}_2)\bd{T}_2 \Big] b\xi_1^2 \xi_2^2 \nonumber\\
	&+ O(|\xi|^6 + b^6).
\end{align}

Substituting \eqref{w_stress} into \eqref{e_stress_1}, and changing to polar coordinates as before, we get
\begin{align}
    	\epsilon =& -\frac{3}{4\pi\del^3} \int_0^{2\pi}\int_0^{\infty} \bd{m}(\del\zeta,\del\lambda) \frac{\phi_3(\sqrt{\eta^2+\lambda^2})}{(\eta^2+\lambda^2)^{5/2}} \eta d\eta d\theta,
\end{align}
which simplifies to
\begin{subequations}
\begin{align}
	\label{e_stress_last_a}
	\epsilon = -\frac{3}{4} \del \Big\{&- [(\bd{q}_i\cdot\bd{T}_i)\bd{n}_0 + (\bd{q}_i\cdot\bd{n}_0)\bd{T}_i] \left(1+\del\lambda(H+\kappa_i)\right)I_{3a} \\
	\label{e_stress_last_b}
	&+ \frac{1}{2}(\bd{q}_{ii}\cdot\bd{n}_0)\bd{n}_0 \del\lambda I_{3a} \\
	\label{e_stress_last_c}
	&+ \frac{1}{8}\Big[ (\bd{q}_{ii}\cdot\bd{T}_j)\bd{T}_j + 2(\bd{q}_{ij}\cdot\bd{T}_i)\bd{T}_j \Big] \del \lambda I_{3b} \Big\} + O(\del^3),
\end{align}
\end{subequations}
with $H$ the mean curvature, as before, and
\begin{subequations}
\begin{align}
	\label{I5}
   	I_{3a} =& \int_0^{\infty} \frac{\phi_3(\sqrt{\eta^2+\lambda^2})}{(\eta^2+\lambda^2)^{5/2}} \lambda^2\eta^3 d\eta, \\
    	\label{I6}
    	I_{3b} =& \int_0^{\infty} \frac{\phi_3(\sqrt{\eta^2+\lambda^2})}{(\eta^2+\lambda^2)^{5/2}} \eta^5 d\eta.
\end{align}
\end{subequations}

The quantities in the square brackets in \eqref{e_stress_last_a} are the surface divergence times the normal, $(\nabla_S\cdot\bd{q})\bd{n}_0$, and 
the surface gradient $\nabla_S(\bd{q}\cdot\bd{n})$.  These are similar to
the Stokeslet terms except now we have $\bd{q}_0=0$ because of the subtraction. 
We will show that the sum of terms with $\kappa_i$ in \eqref{e_stress_last_a} and
\eqref{e_stress_last_b} is $\delta\lambda I_{3a}/2$ times
$\Delta_S[(\bd{q}\cdot\bd{n})\bd{n}]$. We compute the latter using the facts that
$\bd{q}_0 = 0$,  and in our special coordinates $\Delta_S = \Sigma_i \pa_i^2$ and
$\pa_i\bd{n} = -\kappa_i\bd{T}_i$, both at $0$.  We find
\beq \pa_i[(\bd{q}\cdot\bd{n})\bd{n}] = (\bd{q}_i\cdot\bd{n})\bd{n}
          - \kappa_i(\bd{q}\cdot\bd{T}_i)\bd{n}           
          - \kappa_i(\bd{q}\cdot\bd{n})\bd{T}_i  \eeq
and
\beq \pa_i^2[(\bd{q}\cdot\bd{n})\bd{n}] = (\bd{q}_{ii}\cdot\bd{n})\bd{n}
          - 2\kappa_i(\bd{q}_i\cdot\bd{T}_i)\bd{n}           
          - 2\kappa_i(\bd{q}_i\cdot\bd{n})\bd{T}_i  \eeq
The latter matches \eqref{e_stress_last_a}-\eqref{e_stress_last_b}, verifying our assertion.  We note for later use that
\beq 
\label{lapl_n_tan}
	\left\{ \Delta_S[(\bd{q}\cdot\bd{n})\bd{n}] \right\}^{{\rm tan}} = 
        - 2\kappa_i(\bd{q}_i\cdot\bd{n})\bd{T}_i  
\eeq
with sum over $i$, where ${\rm tan}$ means the tangential part.

Since $\Delta_S \bd{q} = \Sigma_i \bd{q}_{ii}$, the first term in brackets in \eqref{e_stress_last_c} is $\{\Delta_S \bd{q}\}^{{\rm tan}}$. We will relate the second term to the
surface gradient of the surface divergence of $\bd{q}$. To identify the latter,
we first note that
\beq \label{div_q} \nabla_S\cdot \bd{q} = \Sigma_i\,\pa_i\bd{q}^{{\rm tan}} \cdot \bd{T}^*_i = \Sigma_i\,\pa_i\bd{q} \cdot \bd{T}^*_i
     -  \Sigma_i\,\pa_i\bd{q}^{{\rm nl}} \cdot \bd{T}^*_i
\eeq 
where
$\bd{q}^{{\rm nl}} = (\bd{q}\cdot\bd{n})\bd{n}$ is the normal part of $\bd{q}$.
Also, for scalar $F$ we have $\nabla_S F = \Sigma_j (\pa_j F) \bd{T}^*_j$.
To find $\nabla_S( \nabla_S\cdot \bd{q} )$ at $0$, we begin with
\beq \nabla_S(\bd{q}_i\cdot \bd{T}^*_i ) = [\bd{q}_{ij}\cdot\bd{T}_i + 
          \bd{q}_i\cdot(\kappa_i\bd{n})\delta_{ij}]\bd{T}_j
   = (\bd{q}_{ij}\cdot\bd{T}_i)\bd{T}_j + \kappa_i(\bd{q}_i\cdot\bd{n})\bd{T}_i
\eeq
with $\bd{q}_i = \pa_i\bd{q}$ etc. and sums over $i,j$ implied.  Then for the normal part we have
\begin{multline}
\label{grad_nl_dot_T}
 \nabla_S(\pa_i\bd{q}^{{\rm nl}}\cdot \bd{T}^*_i ) = 
       \nabla_S[ (\bd{q}\cdot\bd{n})(\pa_i\bd{n})\cdot \bd{T}^*_i ] \\
       = \nabla_S[ (\bd{q}\cdot\bd{n})(-\kappa_i\bd{T}_i)\cdot \bd{T}^*_i ]
       = - 2H \nabla_S(\bd{q}\cdot\bd{n})
\end{multline}
where in the first step we have used the fact that $\bd{n}\cdot\bd{T}^*_i \equiv 0$. Now combining \eqref{div_q}-\eqref{grad_nl_dot_T} and using \eqref{lapl_n_tan} we get
\beq \nabla_S( \nabla_S\cdot \bd{q} ) = (\bd{q}_{ij}\cdot\bd{T}_i)\bd{T}_j
 -\frac12  \left\{ \Delta_S \bd{q}^{{\rm nl}} \right\}^{{\rm tan}}
    +  2H \nabla_S(\bd{q}\cdot\bd{n})\eeq
In summary we have shown that the terms in brackets in \eqref{e_stress_last_c} equal
\beq \{\Delta_S \bd{q}\}^{{\rm tan}}
       + \left\{\Delta_S \bd{q}^{{\rm nl}}\right\}^{{\rm tan}}
       + 2\nabla_S( \nabla_S\cdot \bd{q} ) 
      - 4 H \nabla_S(\bd{q}\cdot\bd{n})\eeq
We can now write a complete formula for $\epsilon$, leading to the correction \eqref{Corr_DL},
\begin{multline}
	\epsilon = \frac34 \del (1+\lambda\del H)I_{3a}
         \Big[ \nabla_S(\bd{q}\cdot\bd{n}) + (\nabla_S\cdot\bd{q})\bd{n}_0 \Big]
         - \frac38 \del^2\lambda I_{3a} [\Delta_S \bd{q}^{{\rm nl}}]    \\
         - \frac{3}{32}\del^2\lambda I_{3b} \Big[
               \{ {\Delta_S (\bd{q}} + \bd{q}^{{\rm nl}}) \}^{{\rm tan}}
           + 2\nabla_S( \nabla_S\cdot \bd{q} ) - 4H \nabla_S(\bd{q}\cdot\bd{n})\Big] + O(\del^3).
\end{multline}


\subsection{Evaluation on the boundary}

When evaluating the integrals on the boundary, we can remove the lower order terms in the regularization error and make it $O(\del^5)$ by an appropriate choice of the smoothing function. To find this
function, we write $s^{\#}(r) = s(r) + a r s'(r) + b r^2 s''(r)$, where constants $a$ and $b$ are chosen to make two moments equal to 0. In \eqref{eps_eta}-\eqref{I3} we now have $\lambda=0$. For the first part of the Stokeslet, with $\phi_1$, the moment conditions are $I_1=0$ and a similar integral with $\eta^3$ in place of $\eta$. This gives us the function used in \cite{beale16} for the Laplacian single layer, 
\beq
	\label{s1_high}
	s_1^{\#}(r) = \erf(r) - \frac23 r (2 r^2 - 5) e^{-r^2}/\sqrt{\pi}.
\eeq

For the second part of the Stokeslet (see \eqref{eps_eta}-\eqref{I3} with $\phi_2$), the $I_{2a}$ term vanishes on the surface, and the moment conditions are $I_{2b}=0$ and a similar integral with $\eta^5$. The fifth order smoothing is
\beq
	\label{s2_high}
	s_2^{\#}(r) = \erf(r) - \frac23 r (4 r^4 - 14 r^2 + 3) e^{-r^2}/\sqrt{\pi}.
\eeq
A similar approach for the stresslet leads us to
\beq
	\label{s3_high}
	s_3^{\#}(r) = \erf(r) - \frac{2}{9} r (8 r^6 - 36 r^4 + 6 r^2 + 9) e^{-r^2}/\sqrt{\pi}.
\eeq

Using these regularization functions for the case on the surface gives high order convergence without the need to compute corrections.


\section{Numerical experiments}

\subsection{Monge parameterization}

In the table below we describe some details on how the geometric quantities were computed. Much of this can be found in Appendix B of \cite{beale16}, but we summarize it here for completeness. Suppose the surface is given by $\phi(x_1,x_2,x_3)=0$, with $\phi>0$ outside, and the normal vector $\nabla\phi$ is predominantly in the $x_3$ direction. Then the parameterization is given by $x_3=z(x_1,x_2)$, where $z$ is the vertical coordinate on the surface.\\

\begin{tabular}{l|l|l}
	Quantity & Evaluation formula & Notes\\
	\hline
	$z_i\equiv \pa z/\pa x_i$ & $-\phi_i/\phi_3$ & by implicit differentiation \\[12pt]
	Tangents & $\bd{T}_1 = (1,0,z_1)$ & \\
	& $\bd{T}_2 = (0,1,z_2)$ & \\[12pt]
	Metric & $g_{ij} = \left( \begin{array}{cc} 
	1+z_1^2 & z_1z_2 \\
	z_1z_2 & 1+z_2^2
	\end{array}\right)$ & \\[12pt]
	Determinant & $g=|g_{ij}| = 1+z_1^2+z_2^2$ & \\[12pt]
	Inverse metric & $\displaystyle g^{ij} = \frac{1}{g}\left( \begin{array}{cc} 
	1+z_2^2 & -z_1z_2 \\
	-z_1z_2 & 1+z_1^2
	\end{array}\right)$ & \\[12pt]
	Dual tangents & $\bd{T}_i^* = \sum_{j=1}^2 g^{ij} \bd{T}_j$, $i=1,2$ & \\[12pt]
	Outward normal & $\bd{n} = \nabla\phi / |\nabla\phi|$ or & \\
	& $\bd{n} = \pm(-z_1,-z_2,1)/\sqrt{g}$ & $+$ if $x_3>z(x_1,x_2)$ outside and \\
	& & $-$ otherwise \\[12pt]
	Mean curvature & $2H = -(\phi_{ii}\phi_j^2 - \phi_i\phi_j\phi_{ij})/\lvert \nabla \phi \rvert^3$ & \\[12pt]
	Stokeslet density & $\bd{f}^{\rm nl} = \bd{f}\cdot\bd{n}\bd{n}$, $\bd{f}^{\rm tan} = \bd{f}-\bd{f}^{\rm nl}$ & \\[12pt]
	Stresslet density & $\bd{q}^{\rm nl} = \bd{q}\cdot\bd{n}\bd{n}$, $\bd{q}^{\rm tan} = \bd{q}-\bd{q}^{\rm nl}$ & \\[12pt]
	Surface gradient & $\displaystyle \nabla_S f = f_1 \bd{T}_1^* + f_2 \bd{T}_2^*$ & $f_i\equiv \pa f/\pa x_i$ by interpolation \\[12pt]
	Surface divergence & $\displaystyle \nabla_S\cdot\bd{v} = \bd{v}^{\rm tan}_1 \cdot \bd{T}_1^* + \bd{v}^{\rm tan}_2 \cdot \bd{T}_2^*$ & $\bd{v}^{\rm tan}_i\equiv \pa \bd{v}^{\rm tan}/\pa x_i$ by interp. \\[12pt]
	Surface Laplacian & $\displaystyle \Delta_S f = \sum_{i,j=1}^2 g^{ij} f_{ij} + \sum_{i=1}^2 c_i f_i$ & $f_i$ and $f_{ij}\equiv \pa^2 f/\pa x_i x_j$ by interp. \\
	Coefficient $c_i$ & $\displaystyle c_i = \mp 2H z_i/\sqrt{g}$ & \cite{beale16} Appendix B\\[12pt]
\end{tabular}


\subsection{Flow around a translating spheroid}

We first define
\beq
	\phi(x_1,x_2,x_3) = \frac{x_1^2}{a^2} + \frac{x_2^2+x_3^2}{b^2} - 1,
\eeq
and compute the flow around a sphere ($a=b$) and a prolate spheroid ($a>b$), translating with velocity $U$. This allows us to test the single layer integral alone. The exact solutions for the flow are well known for both the sphere and the spheroid \cite{chwang75,liron-barta}. Taking $U=(1,0,0)$, the velocity outside the spheroid is 
\beq
	\label{vel_exact}
	\bd{u}(\bd{x}) = 2\alpha B_1 \bd{e}_1 + \alpha r \Big( \frac{1}{R_2} - \frac{1}{R_1} \Big) \bd{e}_r - \alpha r^2 B_2\bd{e}_1 + 2\beta \nabla B_3,
\eeq
where
\begin{align}
	&B_1 = \log \frac{R_2-(x_1-c)}{R_1-(x_1+c)}, \qquad B_2 = \frac{1}{r^2}\Big(\frac{x_1+c}{R_1}-\frac{x_1-c}{R_2} \Big), \qquad B_3 = R_2-R_1+x_1 B_1, \\
	&R_1 = \sqrt{(x_1+c)^2 + r^2}, \qquad R_2 = \sqrt{(x_1-c)^2 + r^2}, \qquad r^2 = x_2^2+x_3^2, \\
	&\alpha = \frac{e^2}{D_e},\qquad \beta=\alpha\frac{1-e^2}{2e^2}, \qquad D_e = (1+e^2)\ln\frac{1+e}{1-e} - 2e, 
\end{align}
$e=c/a$, $c^2=a^2-b^2$, $\bd{e}_r = (x_2\bd{e}_2+x_3\bd{e}_3)/r$ is the unit radial vector in the $x_2x_3$-plane, and $\bd{e}_i$, $i=1,2,3$ are the unit basis vectors. The surface traction is given by $(f_1(\bd{x}),0,0)$ with 
\beq
	\label{f}
	f_1(\bd{x}) = -\frac{4 a e^3}{bD_e \sqrt{a^2-(ex_1)^2}}.
\eeq
The expressions for the unit sphere can be found in the references mentioned above, or for example, in \cite{tlupova13}. 

{In this test, we compute the velocity \eqref{SingleLayer} given the surface traction in \eqref{f}. First, we choose a grid size $h$ and find the quadrature points as explained in Section 2. Second, to choose points near the surface we cover $\mathbb{R}^3$ with a 3D grid of size $h$ and select the points that are $\leq h$ distance outside the surface. For these points, we compute the velocity using \eqref{SL_tilde}, \eqref{StokesletHH}, \eqref{Corr_SL}, with $s_1, s_2$ defined in \eqref{s1}, \eqref{s2}, and the regularization parameter $\delta$ chosen so that $\delta/h$ is a constant.  The integral in \eqref{StokesletHH} is discretized by the quadrature rule \eqref{Quadrature}, evaluating the traction $\bd{f}$ at the quadrature points using \eqref{f}.  To compute the velocity on the surface, the same quadrature is used but without corrections and with $s_1, s_2$ in \eqref{StokesletHH} replaced by $s_1^{\#}, s_2^{\#}$ from \eqref{s1_high}, \eqref{s2_high}.  Then we compare the computed velocity with the exact value \eqref{vel_exact}.}  We define the error at a single point as $e(\bd{x}) = \lvert\bd{u}^{\textrm{computed}} (\bd{x}) - \bd{u}^{\textrm{exact}} (\bd{x})\rvert$, where $\lvert\cdot \rvert$ is the vector's Euclidean norm. We then compute either the max or the $L_2$ norm of this error over the evaluation points. The $L_2$ norm is defined as $\lVert e\rVert_2 = \big(\sum_\bd{x} e^2(\bd{x})/n\big)^{1/2}$, where $n$ is the number of evaluation points. All errors reported are absolute errors; the largest velocity in magnitude is the translational velocity $U$, which was taken to be 1 for both the sphere and the spheroid.

\begin{figure}[htb]
\centering
\includegraphics[scale=0.4]{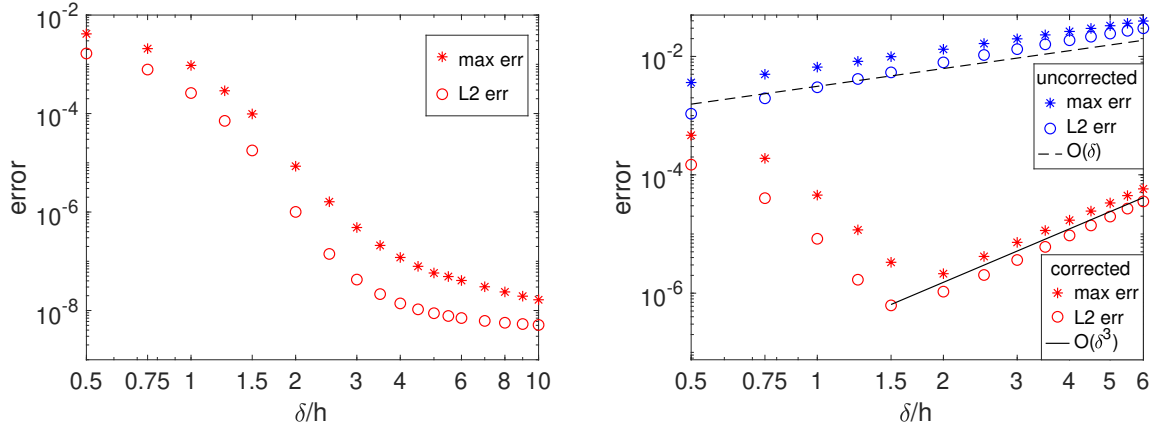}
\includegraphics[scale=0.4]{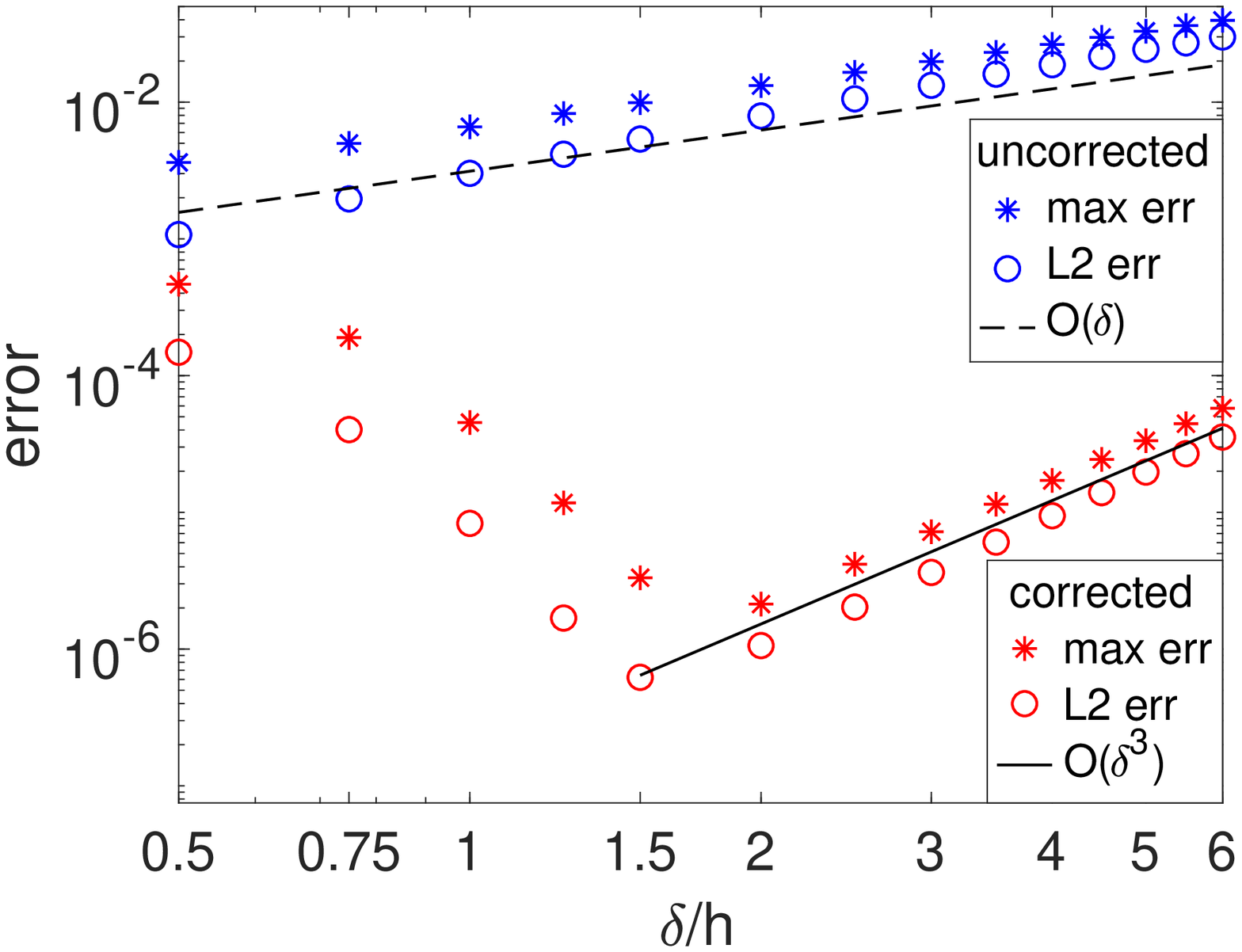}
\caption{Single layer integral, errors for the unit sphere: over quadrature points (left) and over points distance $\leq h$ outside the surface (right).}
\label{SL-sphere-delta}
\end{figure}

{To demonstrate how $\delta/h$ should be chosen in practice, in our first test we keep $h$ fixed and vary $\delta$.} For larger values of $\delta/h$, the discretization error is small and the regularization error is dominant. For small $\delta/h$, the discretization error will dominate, and is expected to be $O(h)$. Figure~\ref{SL-sphere-delta} shows the errors for the unit sphere, where the grid size $h=1/64$ was chosen and $\delta$ was varied. The left graph is the error at the quadrature points, where the higher regularization \eqref{s1_high}-\eqref{s2_high} was used and no corrections are necessary. For a spheroid, rather than a sphere, we find that the error increases with increasing $\delta/h$. Since the integral on the surface is computed with very high accuracy in $\delta$ (Section 3.3), the regularization parameter should be chosen larger to preserve the high order of convergence. The value $\delta=3h$ seems to suffice in all tests we performed. The graph on the right in Figure~\ref{SL-sphere-delta} is the error for points close to the surface. The graph demonstrates a typical behavior where $\delta/h>1$ should be chosen so that the regularization error dominates leading to a higher order of convergence in grid size $h$. The regularization $\delta=2h$ seems sufficiently large to achieve $O(h^3)$ convergence, as shown in our subsequent testing.

\begin{figure}[htb]
\centering
\includegraphics[scale=0.4]{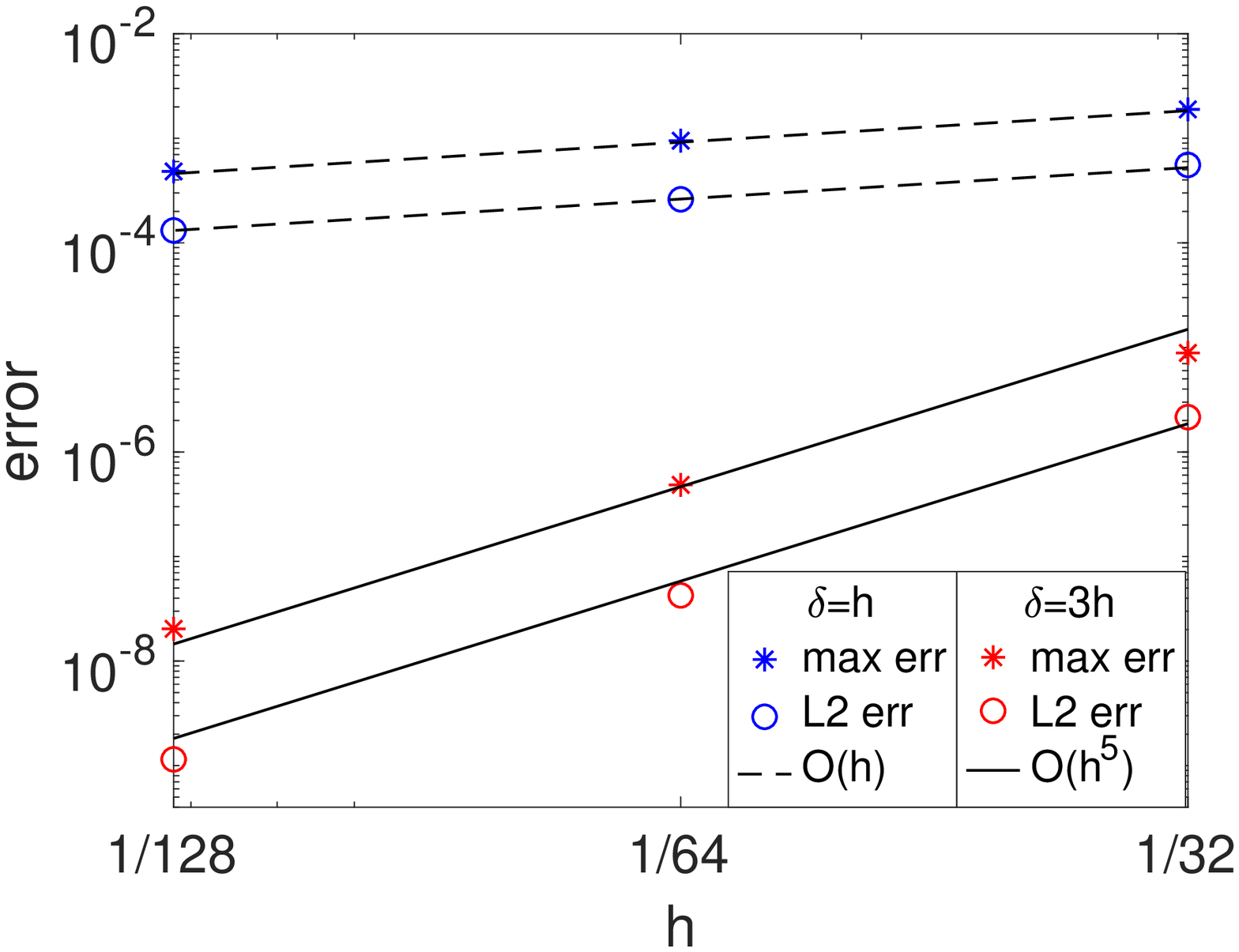}
\includegraphics[scale=0.4]{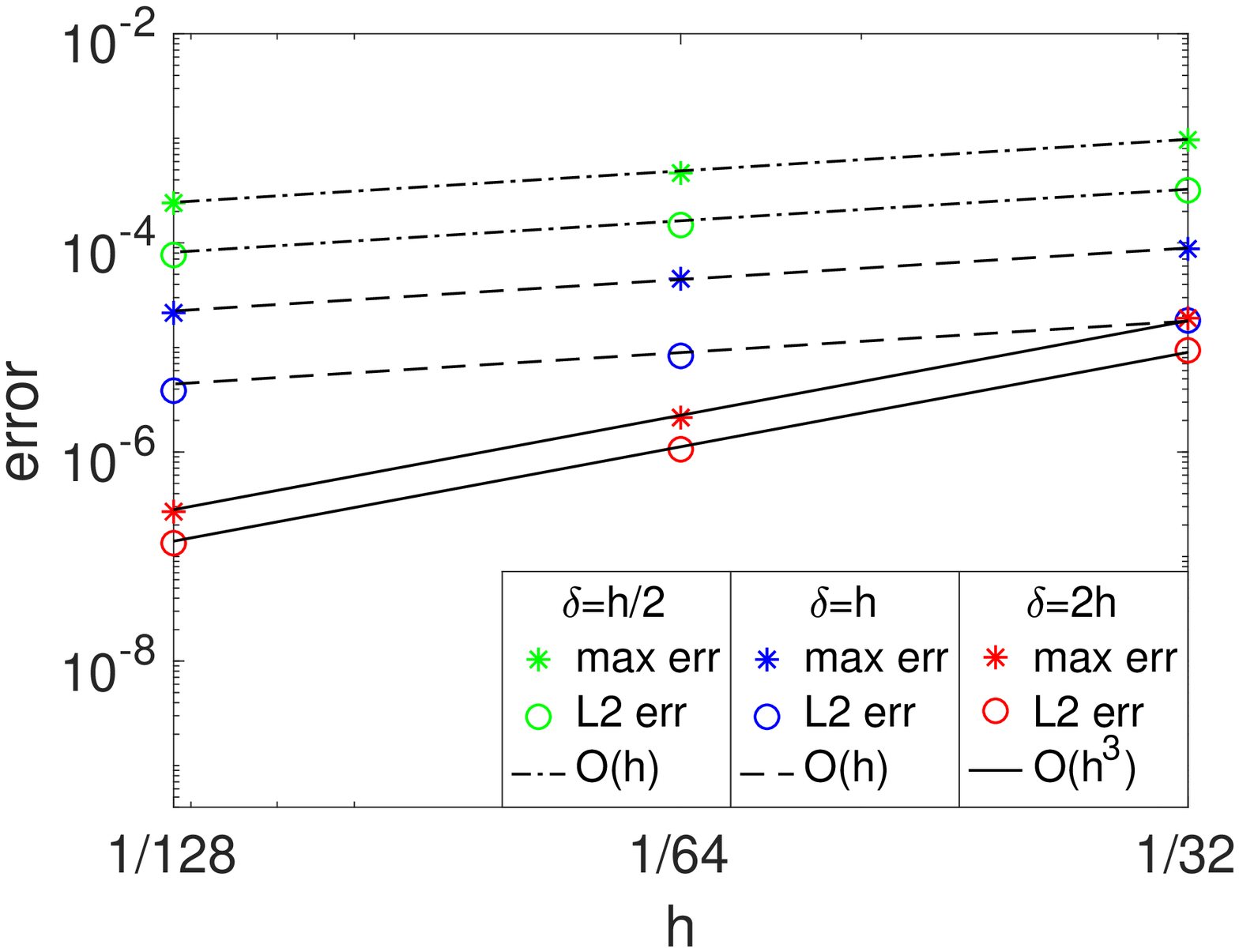}
\caption{Single layer integral, errors for the unit sphere: over quadrature points (left) and over points distance $\leq h$ outside the surface (right).}
\label{SL-sphere}
\end{figure}

\begin{figure}[htb]
\centering
\includegraphics[scale=0.4]{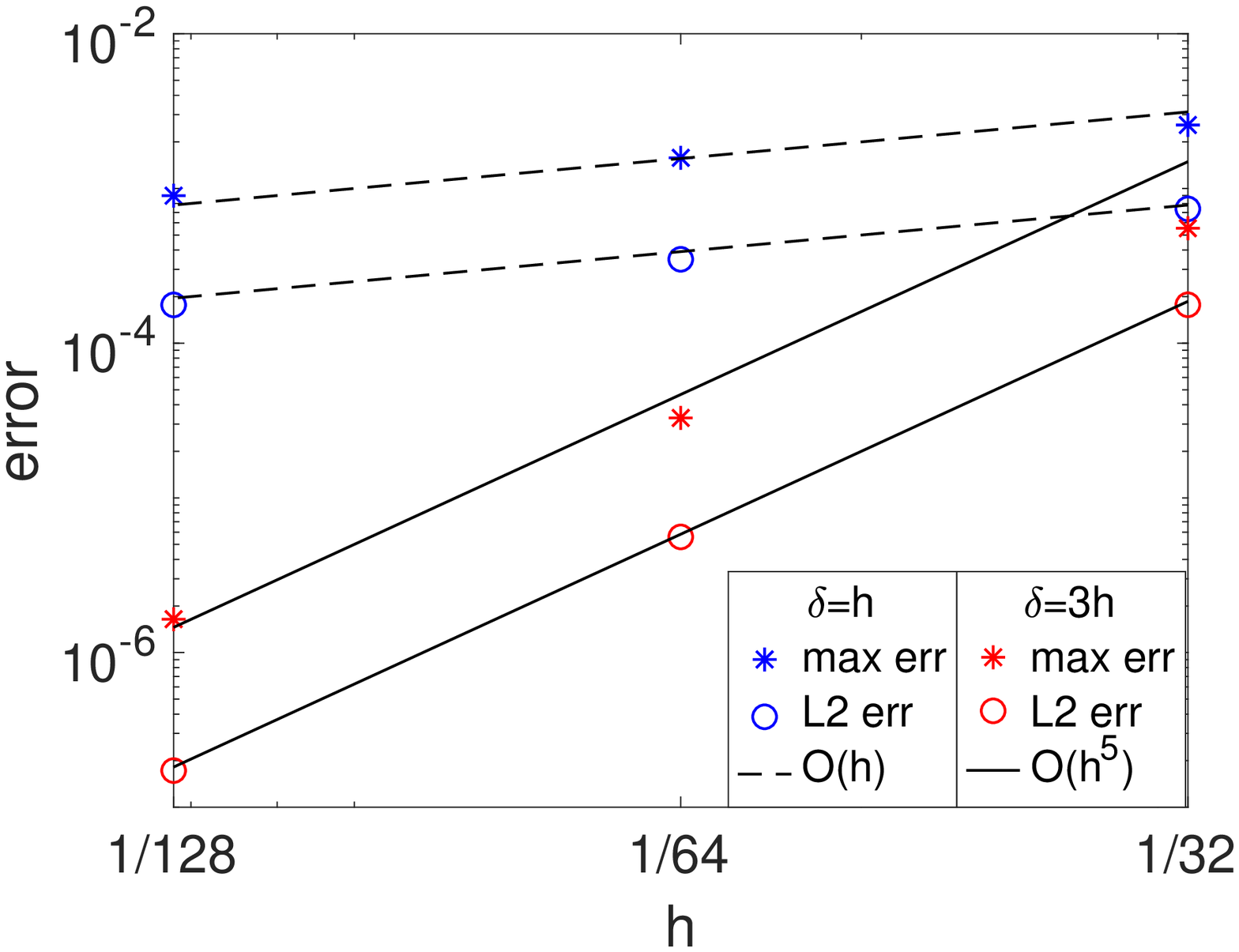}
\includegraphics[scale=0.4]{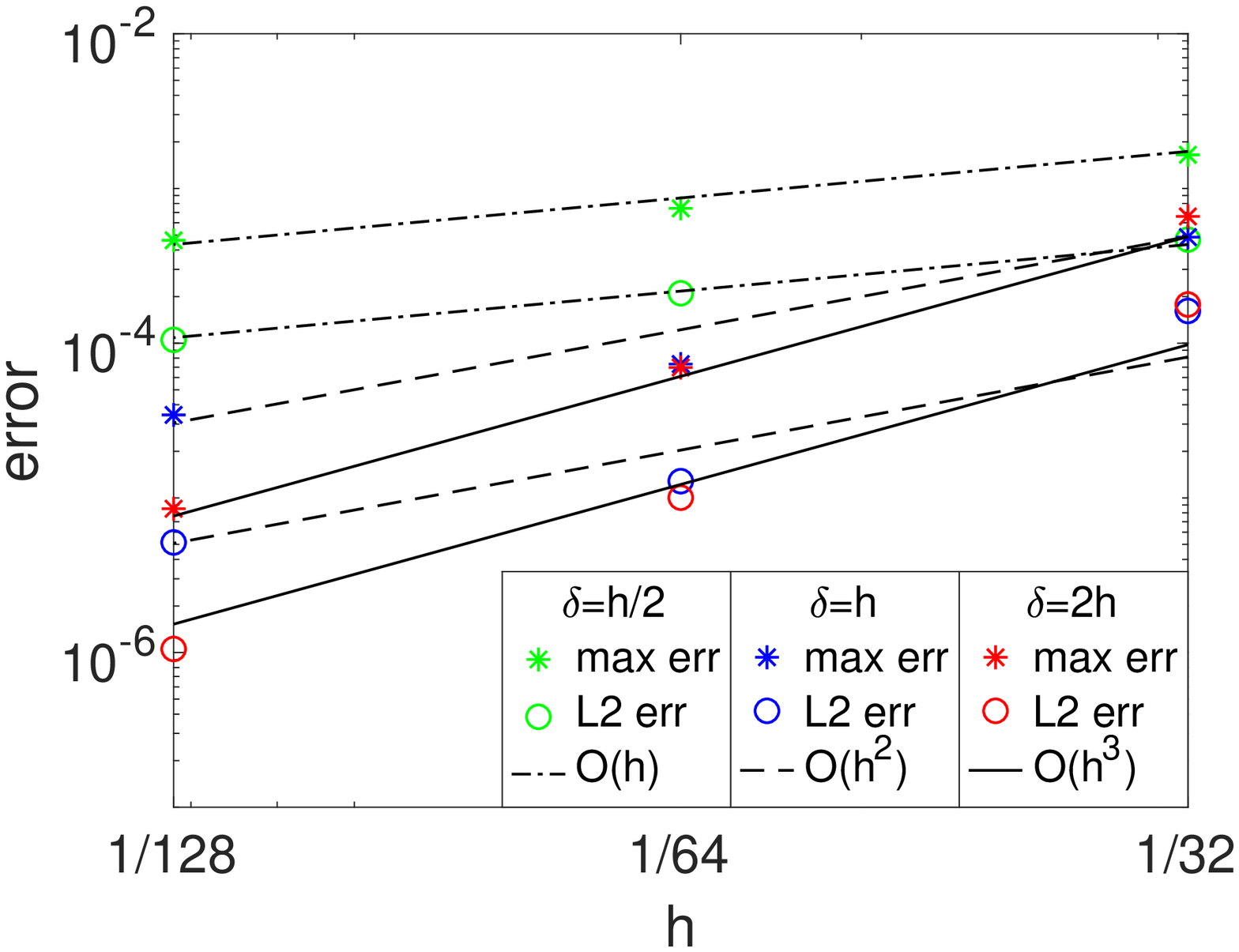}
\caption{Single layer integral, errors for the spheroid $a=1,b=0.5$: over quadrature points (left) and over points distance $\leq h$ outside the surface (right).}
\label{SL-1p5p5}
\end{figure}

We further demonstrate the predicted convergence rates by refining the grid while keeping $\delta/h$ fixed. Figure~\ref{SL-sphere} shows the errors for the unit sphere for points on the surface (taken to be the quadrature points), and points that are $\leq h$ distance outside the surface. We repeat this test for the spheroid $a=1$, $b=0.5$ in Figure~\ref{SL-1p5p5}. The errors are larger for the spheroid due to the larger curvature and varied spacing: for $h=1/32$ for example, there are 17070 quadrature points on the unit sphere, but only 6958 quadrature points on the spheroid. Both cases display the predicted order of convergence. As described earlier, choosing $\delta/h<1$ will result in slower convergence of $O(h)$, while $\delta/h=1$ is in the intermediate regime where convergence can be observed as $O(h)$ (as we see for the sphere) or faster (as it happens to be for the spheroid). A larger $\delta/h=2$ will typically result in the predicted $O(h^3)$ convergence off the surface, as seen for both the sphere and the spheroid. On the surface, an even larger regularization is recommended, and we choose $\delta/h=3$. The errors decrease rapidly with observed convergence rates close to $O(h^5)$.

\subsection{Double layer identity}

\begin{figure}[htb]
\centering
\includegraphics[scale=0.4]{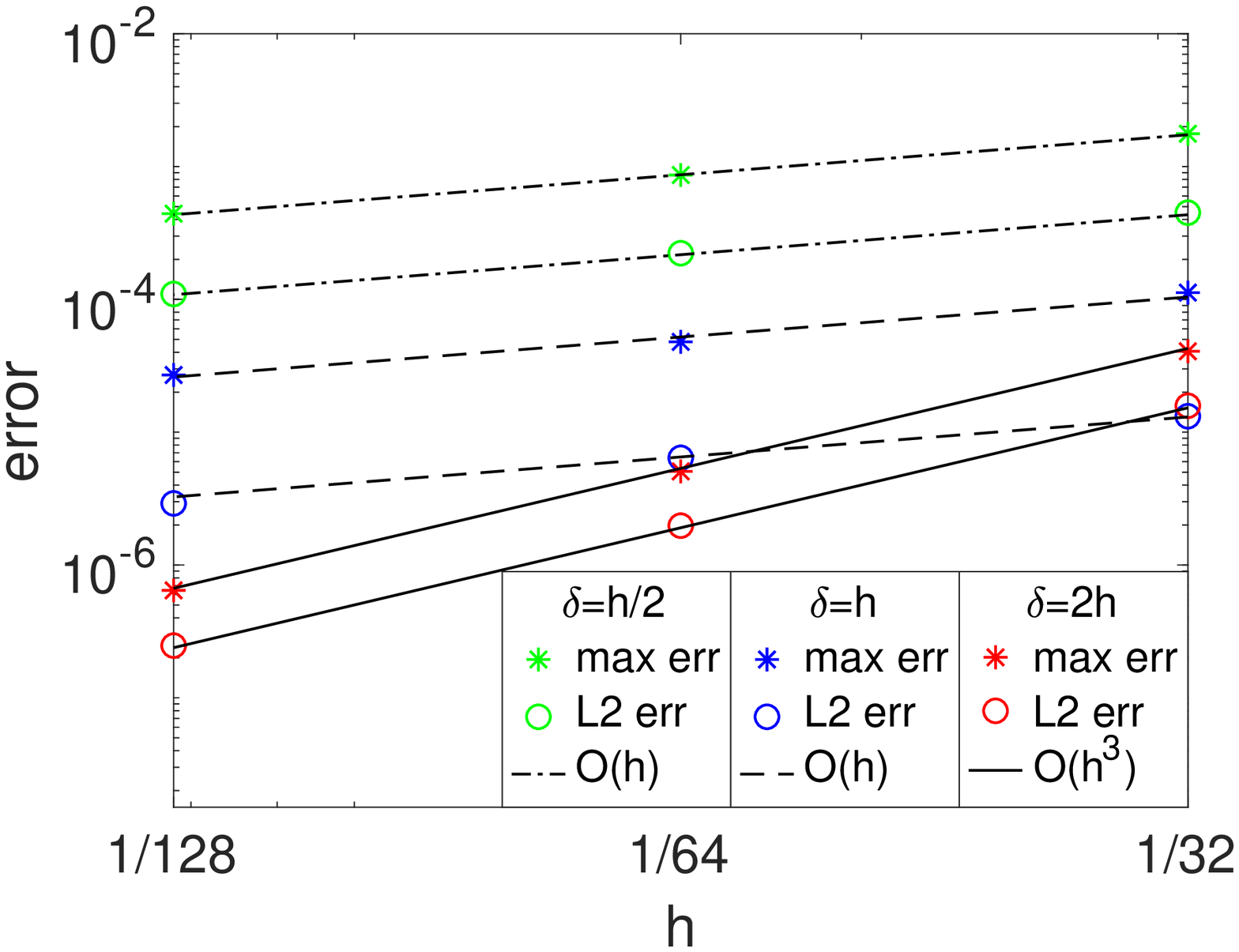}
\includegraphics[scale=0.4]{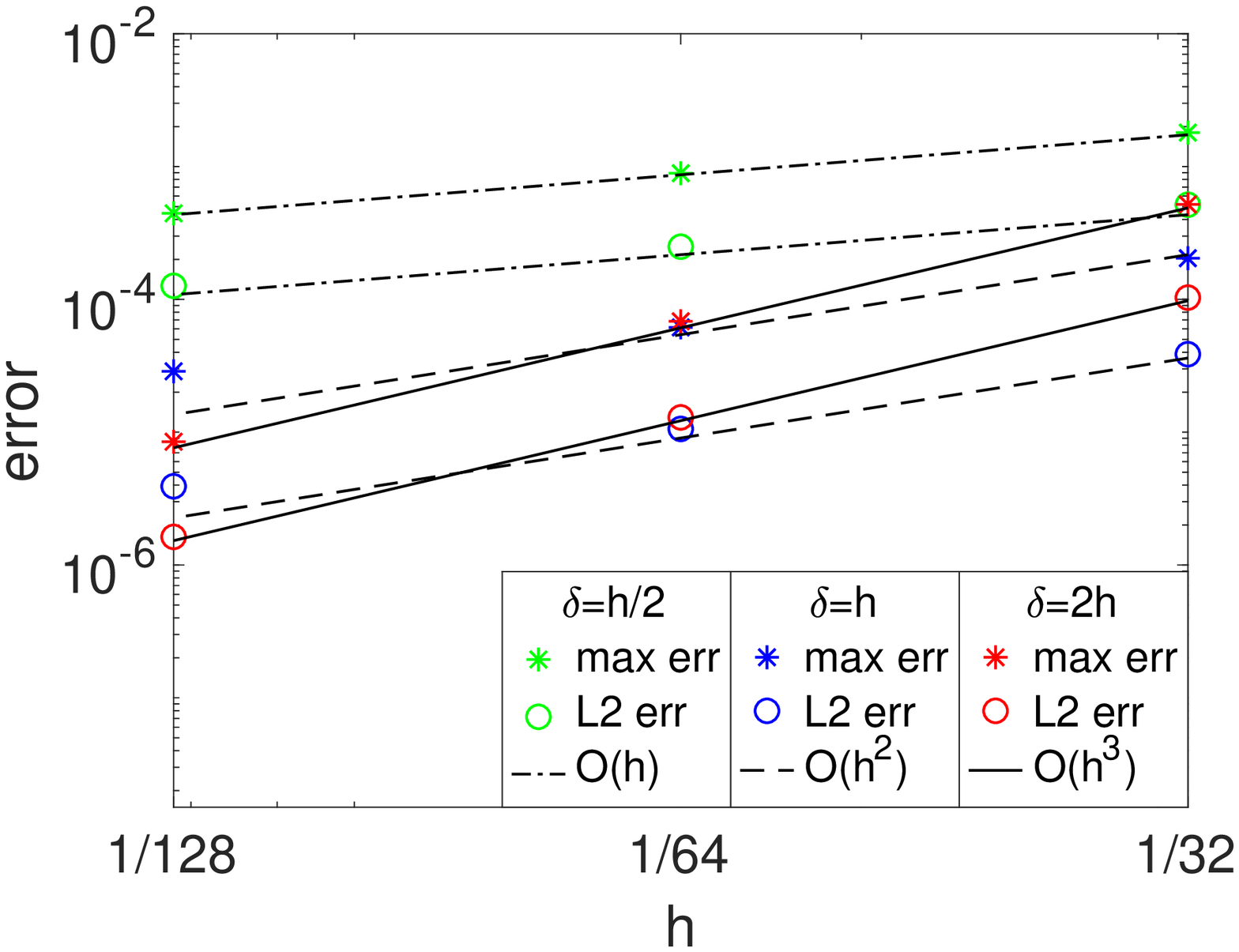}
\caption{Double layer identity, errors over points distance $\leq h$ from the surface for the unit sphere (left) and the ellipsoid $a=1,b=0.6,c=0.4$ (right).}
\label{DL-sphere-1p6p4}
\end{figure}

In order to now test the double layer integral alone, we use the identity (2.3.19) from \cite{pozrikidis92},
\beq 
	\label{DL_identity}
	\frac{1}{8\pi} \epsilon_{ilm} \int_{\partial\Omega} x_m T_{ijk} (\bd{x_0,x}) n_k(\bd{x})dS(\bd{x}) = \chi (\bd{x}_0) \epsilon_{jlm} x_{0,m},
\eeq
where $\chi$ = 1, 1/2, 0 when $\bd{x}_0$ is inside, on, and outside the boundary. We let $l=1$ and define $q_i(\bd{x}) = \epsilon_{i1m}x_m = (0,-x_3,x_2)$. {Given this density function, we compute the left hand side of the identity \eqref{DL_identity} as the double layer potential using \eqref{DoubleLayer_subtract}, \eqref{DL_tilde}, \eqref{StressletHH}, \eqref{Corr_DL}.  We then compare the computed values to the exact values, given by the right hand side of \eqref{DL_identity}.}  It is worth noting that since we use subtraction in the double layer, for this example we get $(\bd{q}-\bd{q}_0)\cdot (\bd{x}-\bd{x}_0) \equiv 0$, and therefore the integral is identically zero on the boundary. For a general ellipsoid, we use
\beq
	\phi(x_1,x_2,x_3) = \frac{x_1^2}{a^2} + \frac{x_2^2}{b^2} + \frac{x_3^2}{c^2} - 1,
\eeq
and test the identity at points inside and outside the surface at distance $\leq h$ away. Figure~\ref{DL-sphere-1p6p4} shows the errors for the unit sphere and the ellipsoid $a=1,b=0.6,c=0.4$. Again, we test three values of $\delta/h$: 0.5, 1, and 2. The behavior is similar to the single layer integral near the surface. The errors reported are again absolute errors, with the largest magnitude of the solution about 1 for the sphere and about $0.6$ for the ellipsoid. The comparison of the two cases is also affected by the differing number of points for given $h$. For reference, the number of quadrature points for the sphere and the ellipsoid for each $h$ is given in Table~\ref{Quad_pts}.


\subsection{Sum of single and double layer}

\begin{figure}[htb]
\centering
\includegraphics[scale=0.4]{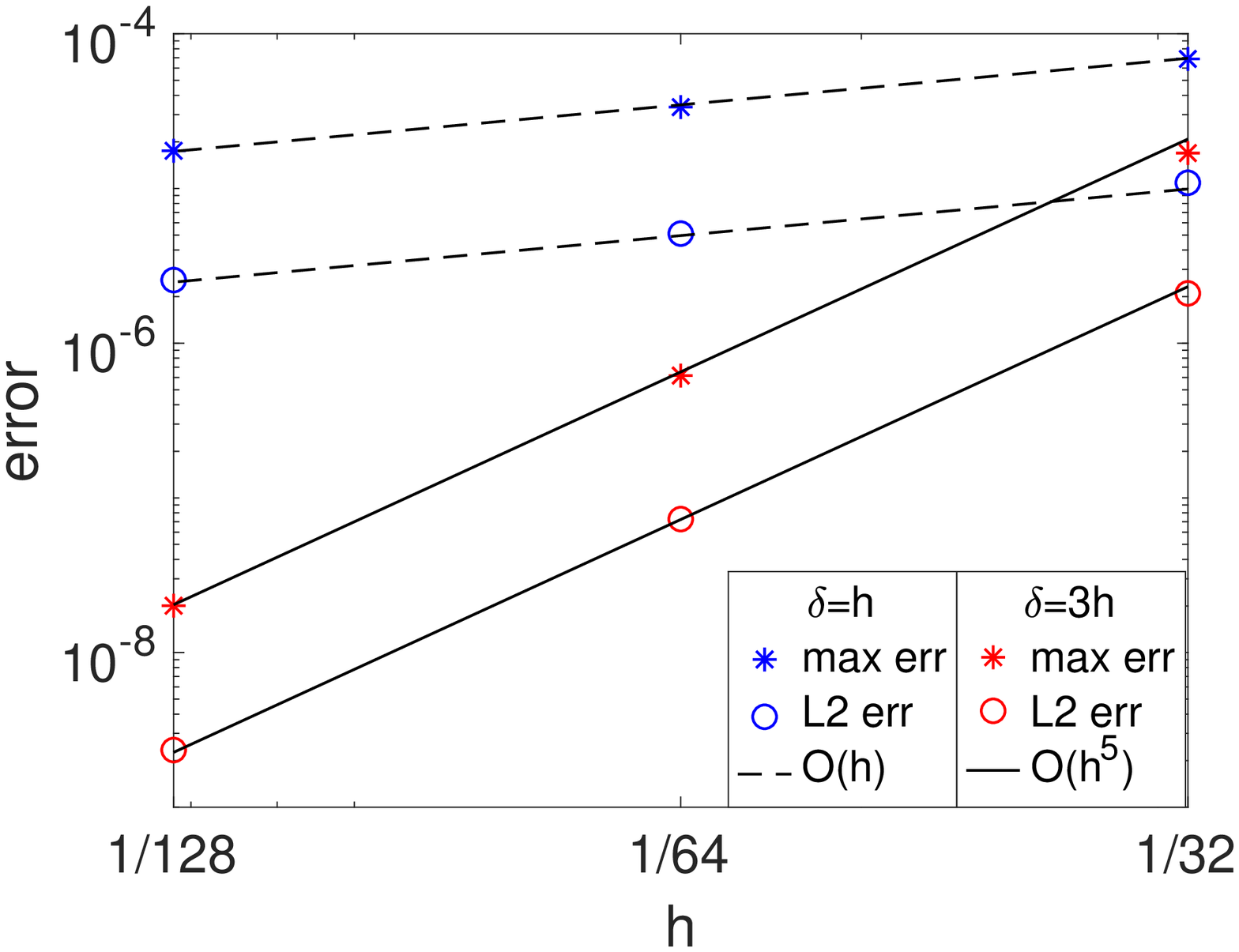}
\includegraphics[scale=0.4]{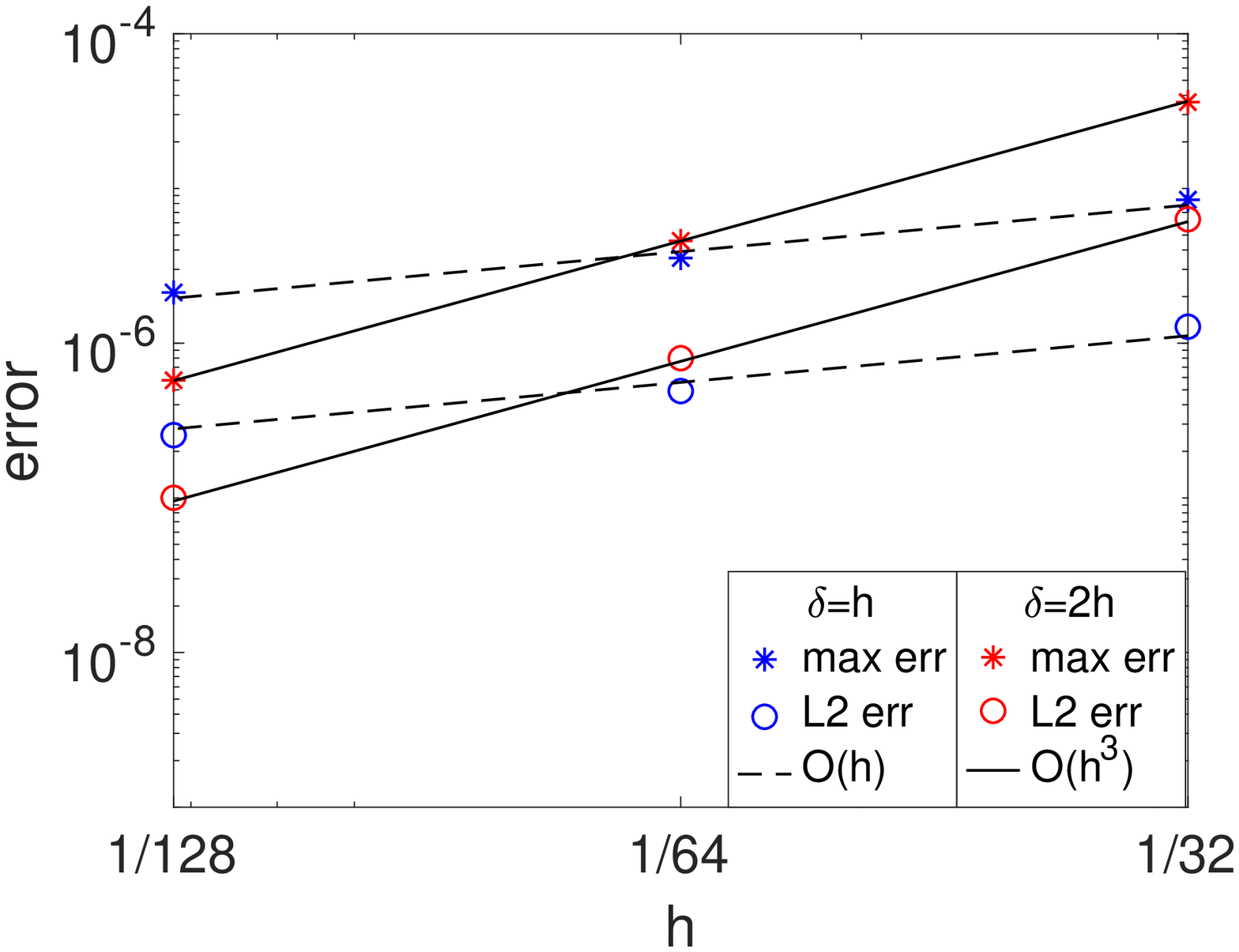}
\caption{Sum of single and double layer, errors for the unit sphere: over quadrature points (left) and over points distance $\leq h$ from the surface (right).}
\label{Sum-sphere}
\end{figure}

\begin{figure}[htb]
\centering
\includegraphics[scale=0.4]{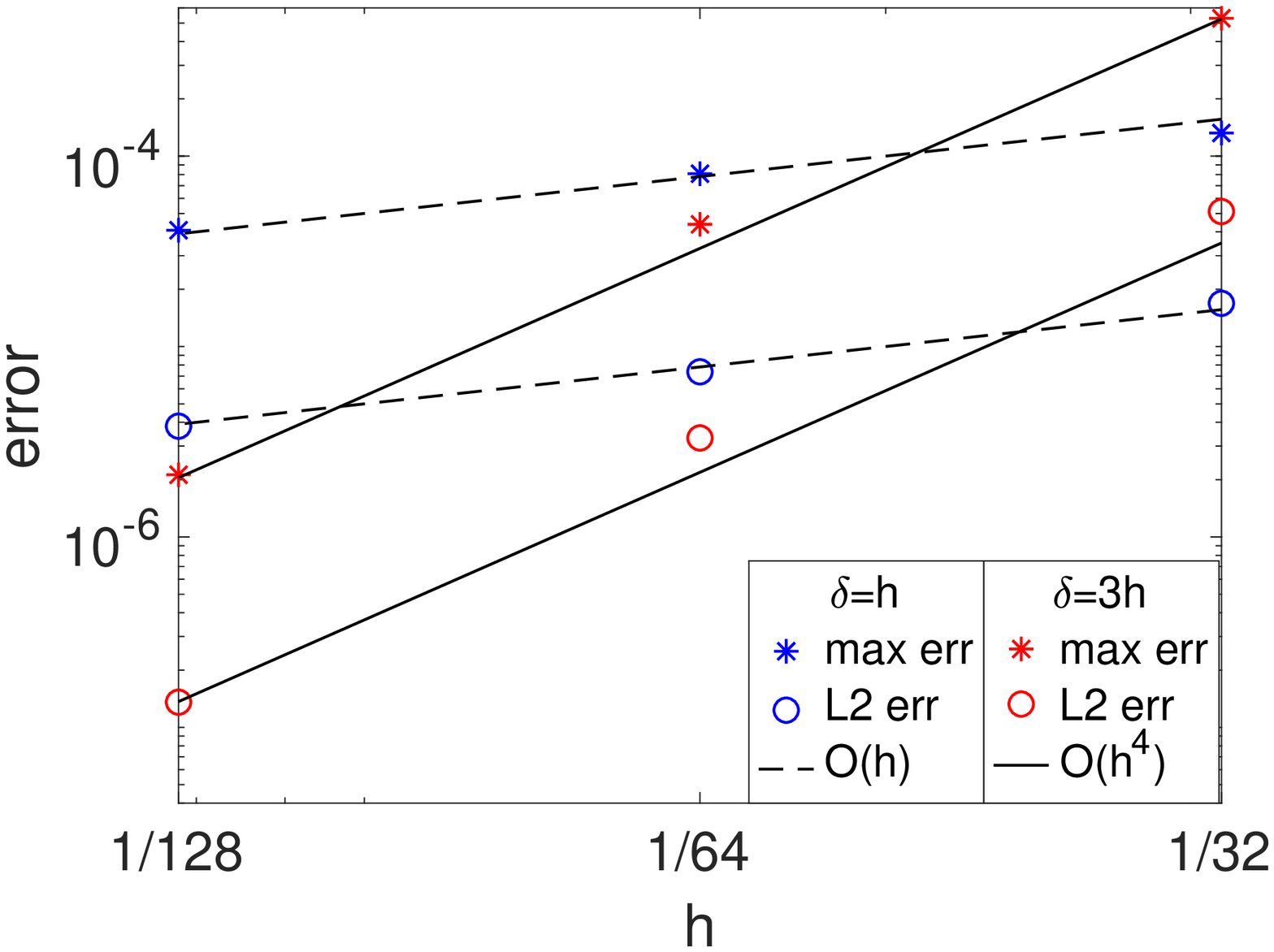}
\includegraphics[scale=0.4]{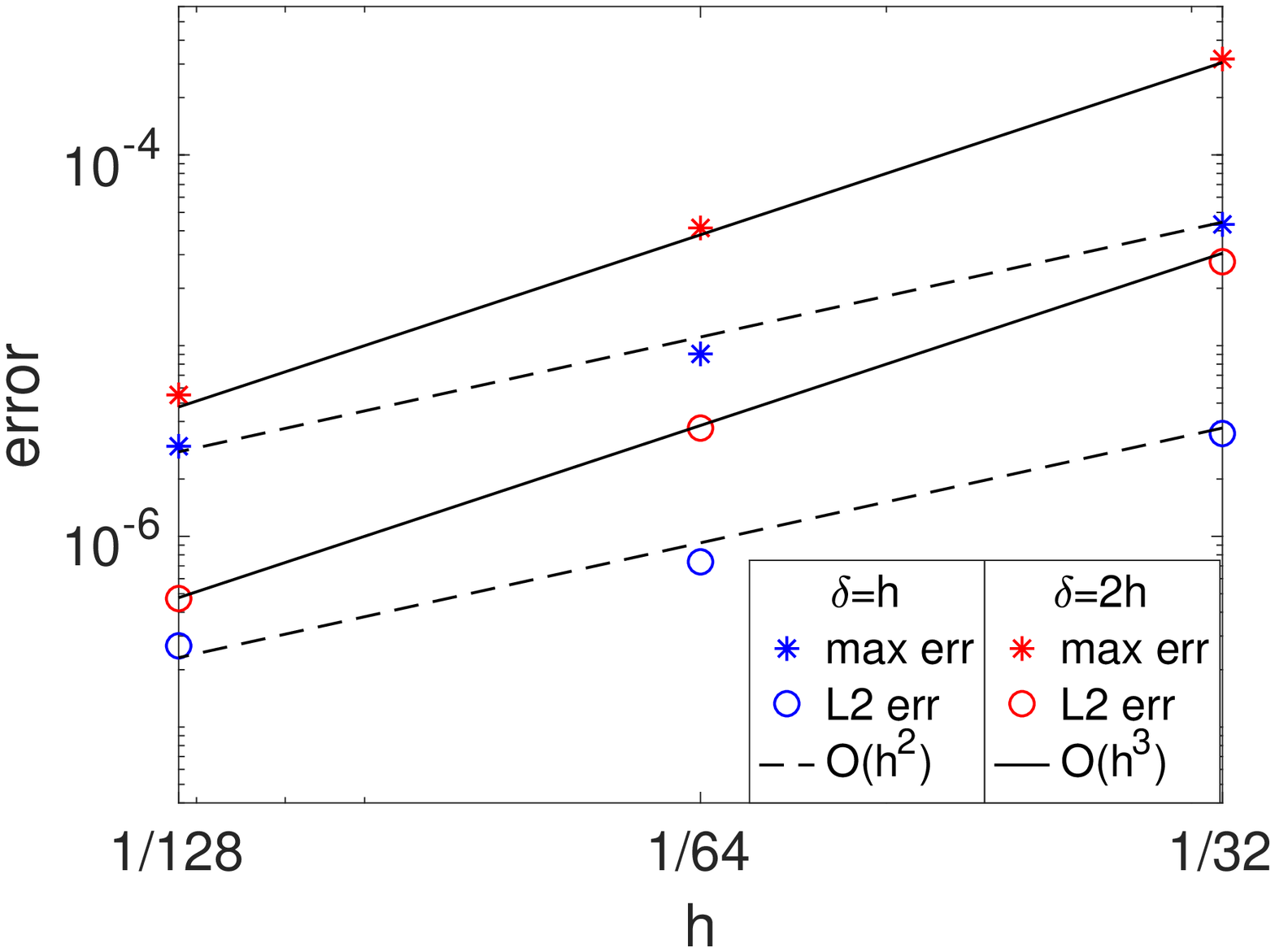}
\caption{Sum of single and double layer, errors for the ellipsoid $a=1,b=0.6,c=0.4$: over quadrature points (left) and over points distance $\leq h$ from the surface (right).}
\label{Sum-1p6p4}
\end{figure}

\begin{figure}[htb]
\centering
\includegraphics[scale=0.4]{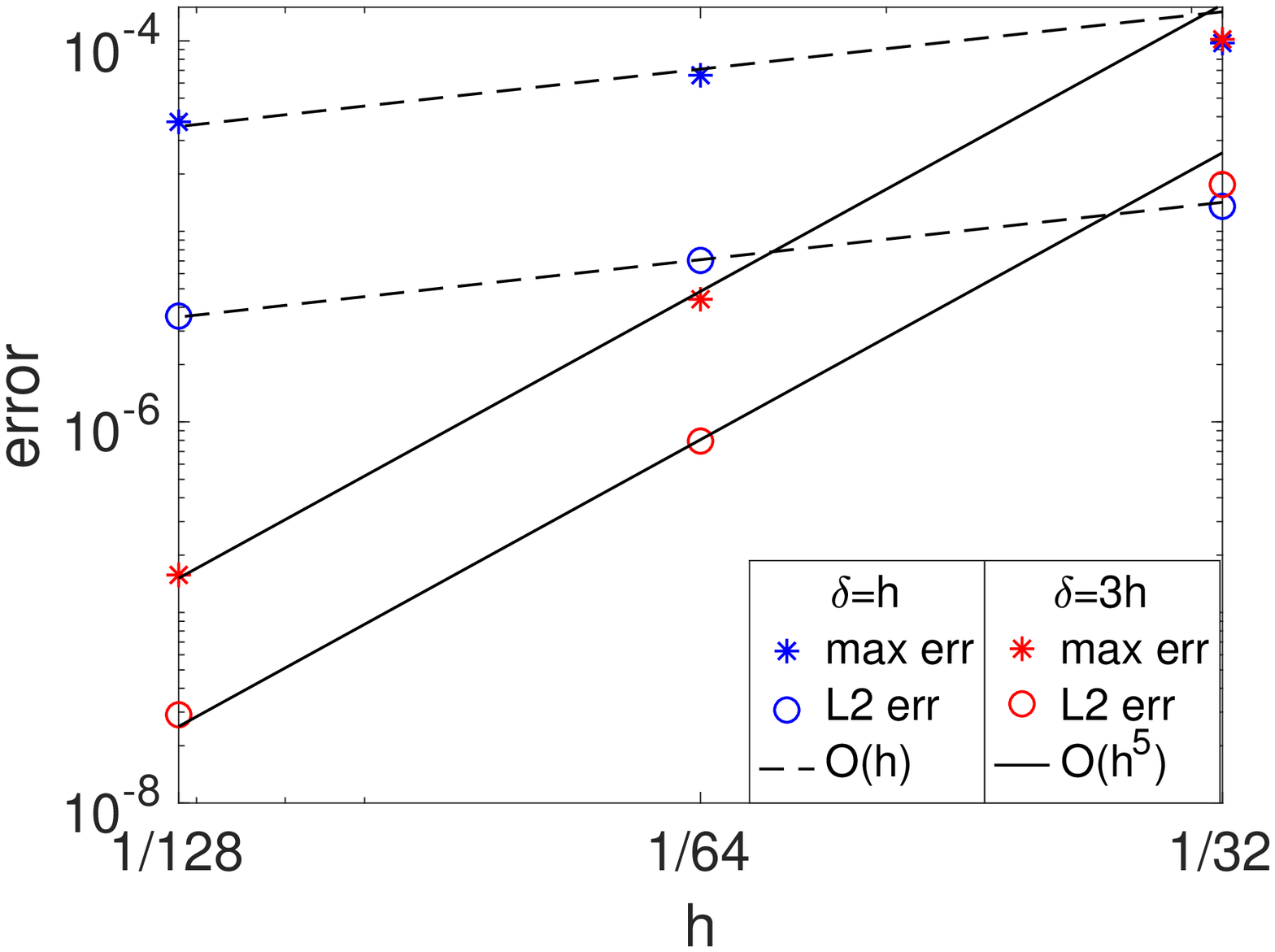}
\includegraphics[scale=0.4]{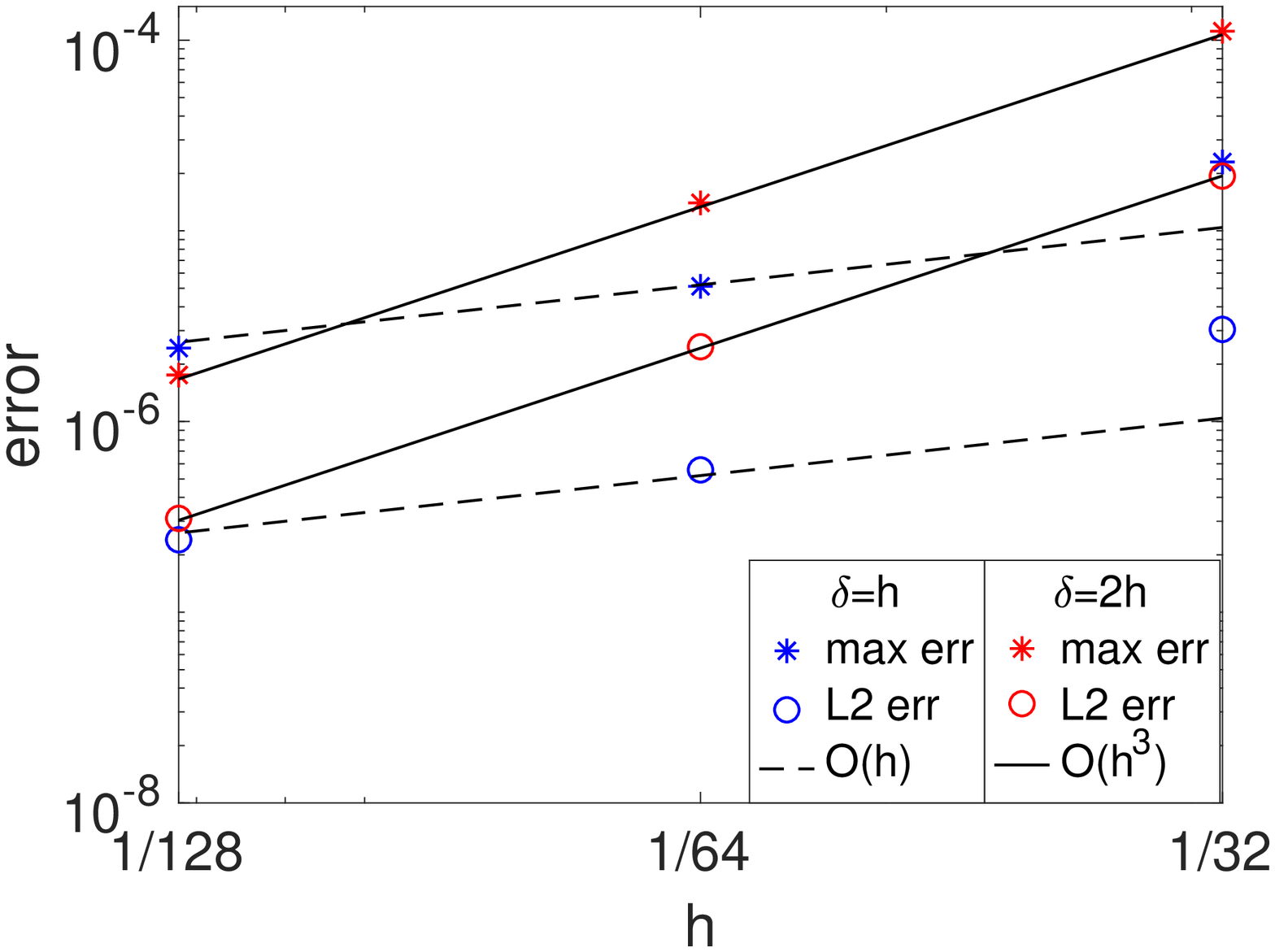}
\caption{Sum of single and double layer, errors for the molecular surface: over quadrature points (left) and over points distance $\leq h$ from the surface (right).}
\label{Sum-molecule}
\end{figure}

One of the advantages of using boundary integral formulations is that jumps in the physical quantities across interfaces get incorporated into the integrals naturally. Specifically, the general integral formulation, expressed as the sum of the single and double layer integrals,
\beq 
	\label{Sum_SLDL}
	u_i(\bd{y}) = 
	-\frac{1}{8\pi}\int_{\partial\Omega} S_{ij} (\bd{y,x}) [f]_j(\bd{x}) dS(\bd{x}) 
	- \frac{1}{8\pi}\int_{\partial\Omega} T_{ijk} (\bd{y,x}) [u]_j(\bd{x}) n_k(\bd{x})dS(\bd{x}),
\eeq
has $[f] = f^+ - f^- = (\sigma^+ - \sigma^-)\cdot \bd{n} $ as the jump in surface force and $[u]$ as the jump in velocity. Here $\bd{n}$ is the outward unit normal, and the plus/minus signs denote the outside/inside of the boundary. To demonstrate the accuracy of this formulation, we use the following solution. On the inside, we assume the velocity is given by a point force singularity of strength $\bd{b} = (1,0,0)$, placed at $\bd{y}_0 = (2,0,0)$. The solution is given by the Stokeslet velocity
\beq 
	\label{Stokeslet_point}
	u^-_i(\bd{y}) = \frac{1}{8\pi} S_{ij}b_j = \frac{1}{8\pi} \Big( \frac{\delta_{ij}}{r} + \frac{\hat{y}_i \hat{y}_j}{r^3} \Big)b_j,
\eeq 
and the stress tensor is 
\beq
	\label{Stress_point}
	\sigma^-_{ik}(\bd{y}) = \frac{1}{8\pi} T_{ijk} b_j = \frac{-6}{8\pi}\frac{\hat{y}_i \hat{y}_j \hat{y}_k}{r^5} b_j,
\eeq
where $\hat{\bd{y}} = \bd{y}-\bd{y}_0$, $r=|\hat{\bd{y}}|$. {We assume this data for the inside of the boundary, and take the solution to be $u^+=0$, $\sigma^+=0$ for the outside.  The jumps $[u]$ and $[f]$ are evaluated at the quadrature points using these inside/outside values. The single layer and double layer integrals in \eqref{Sum_SLDL} are then computed as described in Sections 4.2 and 4.3. To compare to the exact solution on the boundary, we take the solution as the average of outside and inside, or half of the formula for $u_i$ in \eqref{Stokeslet_point}.} This test allows us to check convergence when the formulation involves both the single and double layer potentials, using the high order regularization \eqref{s1_high}-\eqref{s3_high} for points on the surface, and the effect of corrections \eqref{Corr_SL}-\eqref{Corr_DL} for points off the surface. Figure~\ref{Sum-sphere} shows errors for the unit sphere, evaluated at the quadrature points only (left graph) and points inside and outside the surface that are distance $\leq h$ away. Figure~\ref{Sum-1p6p4} shows similar errors for the ellipsoid $a=1,b=0.6,c=0.4$. As another test, we use the four-atom molecular surface as in \cite{beale16}, given by $\sum_{k=1}^4 \exp(- |{\bf x} - {\bf x}_k|^2/r^2) = c$, with centers $(\sqrt{3}/3,0,-\sqrt{6}/12)$, $(-\sqrt{3}/6,\pm .5,-\sqrt{6}/12)$, $(0,0,\sqrt{6}/4)$ and $r = .5$, $c = .6$. The results for this surface are shown in Figure~\ref{Sum-molecule} and exhibit a similar behavior. The largest magnitude of the solution is about $0.079$ for both the sphere and the ellipsoid, and about $0.073$ for the molecular surface. See Table~\ref{Quad_pts} for the number of quadrature points for these surfaces for different grid sizes. On the surface, we see high order convergence, $O(h^5)$ for the sphere and $O(h^4)$ for the ellipsoid, when regularization is chosen large enough, such as $\delta/h=3$. For smaller regularization parameter $\delta/h=1$ or less, observed convergence is $O(h)$, so this is not recommended in practice. For points off the surface, the accuracy in the corrected solution is the predicted $O(h^3)$ for regularization $\delta/h=2$.

\subsection{Flow due to an interface with different viscosities}

Here we consider an example of an interface between two fluids with different viscosities, and an integral equation must be solved to find the interface velocity. The interface undergoes a discontinuity in the surface force $[\bd{f}]$, while the velocity across the interface is continuous \cite{pozrikidis92}. The integral equation for the interface velocity is given by
\begin{align}
	(\lambda+1) u_i(\bd{x}_0) = &-\frac{1}{4\pi \mu_0}\int_{\pa\Omega} S_{ij}(\bd{x}_0,\bd{x}) [f]_j(\bd{x})dS(\bd{x}) \nonumber\\
	&+ \frac{\lambda-1}{4\pi}\int_{\pa\Omega} T_{ijk} (\bd{x}_0,\bd{x}) u_j(\bd{x}) n_k(\bd{x})dS(\bd{x})
\label{IntEq-1surf}
\end{align}
for $\bd{x}_0\in \pa\Omega$, where $\mu_0,\mu_1$ are the external and internal fluid viscosities and $\lambda=\mu_1/\mu_0$. The discontinuity in the surface force is given by $[\bd{f}] = 2\gamma H \bd{n} - \nabla_S\gamma$, where $\gamma$ is the surface tension, $H$ is the mean curvature, and $\bd{n}$ is the outward unit normal \cite{pozrikidis92,tornberg18}. In our numerical tests, we set $\mu_0=1, \mu_1=2$, and $\gamma=1+x_1^2$. We solve the integral equation using successive evaluations, i.e.,
\begin{align}
	(\lambda+1) u^{N}_i(\bd{x}_0) = &-\frac{1}{4\pi \mu_0}\int_{\pa\Omega} S_{ij}(\bd{x}_0,\bd{x}) [f]_j(\bd{x})dS(\bd{x}) \nonumber\\
	&+ \frac{\lambda-1}{4\pi}\int_{\pa\Omega} T_{ijk} (\bd{x}_0,\bd{x}) u^{N-1}_j(\bd{x}) n_k(\bd{x})dS(\bd{x}),
	\label{Successive}
\end{align}
for $N=1,2,...$, and $\bd{u}^0 = \bd{0}$. We stop these iterations when the iteration error, defined as 
\beq
	\label{error_iter}
	e^N := \max_{\bd{x}_0} \lvert \bd{u}^{N}-\bd{u}^{N-1} \rvert,
\eeq 
is below a prescribed tolerance, and $\lvert\cdot \rvert$ is the vector's Euclidean norm. We use the higher order regularization derived in Sec. 3.3. Since the exact solution is not known, we check the convergence rates by defining
\beq
	\label{error_h}
	e_h (\bd{x}) = \bd{u}_{h} (\bd{x}) - \bd{u}_{h/2} (\bd{x}),
\eeq
and taking either the max or the $L_2$ norm of this error over the surface points given by $h$, the larger of the two grid sizes used. These errors are shown in Table~\ref{BIE-1Surf} for the unit sphere and the ellipsoid $a=1,b=0.6,c=0.4$, with $\del=3h$ in both cases. In these tests, it took $N=8$ iterations for the iteration error \eqref{error_iter} to reach below $10^{-10}$ for the sphere and around 12 iterations for the ellipsoid. The error in the solution \eqref{error_h} is larger especially for the ellipsoid, likely coming mostly from evaluating the single layer integral with the surface tension density. This can be remedied somewhat by computing the single layer integral with increased resolution before solving the integral equation, since the single layer is a nonhomogeneous term. For example, we solved the integral equation for each of the values of $h$, but in each case computed the Stokeslet integral at the needed points using the finer grid $h=1/256$. As shown in the last section of Table~\ref{BIE-1Surf}, doing this reduces the error by over an order of magnitude. Table~\ref{BIE-1Surf} reports absolute errors.  The largest velocity magnitude is $0.19$ for the sphere and $0.54$ for the ellipsoid. Again there are fewer points on the ellipsoid than on the sphere for given $h$.

\begin{table}[htb]
\centering
\begin{tabular}{c|r|r|r|r}
\hline
$h$ & Sphere & Spheroid & Ellipsoid & Molecule \\
\hline
1/16 & 4302   & 1766 & 1742  & 2392 \\
\hline
1/32 & 17070   & 6958 & 6902  & 9562 \\
\hline
1/64 & 68166   & 27934 & 27566  & 38354 \\
\hline
1/128 & 272718   & 112006 & 110250  & 153399 \\
\hline
\end{tabular}
\caption{Number of quadrature points: unit sphere, spheroid $a=1,b=c=0.5$, ellipsoid $a=1,b=0.6,c=0.4$, and the molecular surface from \cite{beale16}.}
\label{Quad_pts}
\end{table}

\begin{table}[htb]
\centering
\begin{tabular}{c||c|c||c|c||c|c}
\hline
\multicolumn{1}{c||}{} &
\multicolumn{2}{c||}{Sphere} &
\multicolumn{2}{c||}{Ellipsoid} &
\multicolumn{2}{c}{Ellipsoid*} \\
\hline
$h$ & $\lVert e_h\rVert_\infty$ & $\lVert e_h\rVert_2$ & $\lVert e_h\rVert_\infty$ & $\lVert e_h\rVert_2$ & $\lVert e_h\rVert_\infty$ & $\lVert e_h\rVert_2$ \\
\hline
1/16 & 1.21e-04 & 4.57e-05 & 2.95e-02 & 9.99e-03    & 6.93e-03 & 1.85e-03 \\
\hline
1/32 & 6.24e-06 & 1.59e-06 & 7.35e-03 & 1.55e-03    & 6.86e-04 & 1.20e-04 \\
\hline
1/64 & 3.06e-07 & 4.72e-08 & 7.78e-04 & 1.16e-04    & 5.50e-05 & 7.97e-06 \\
\hline
\end{tabular}
\caption{Flow due to an interface, for the unit sphere and the ellipsoid $a=1,b=0.6,c=0.4$. Ellipsoid*: the single layer integral was computed using $h=1/256$. Grid size $h$, max and $L_2$ norms of the error defined in \eqref{error_h}. Regularization parameter $\del=3h$.}
\label{BIE-1Surf}
\end{table}


\subsection{Two interfaces close to each other}

\begin{figure}[htb]
\centering
\includegraphics[scale=0.5]{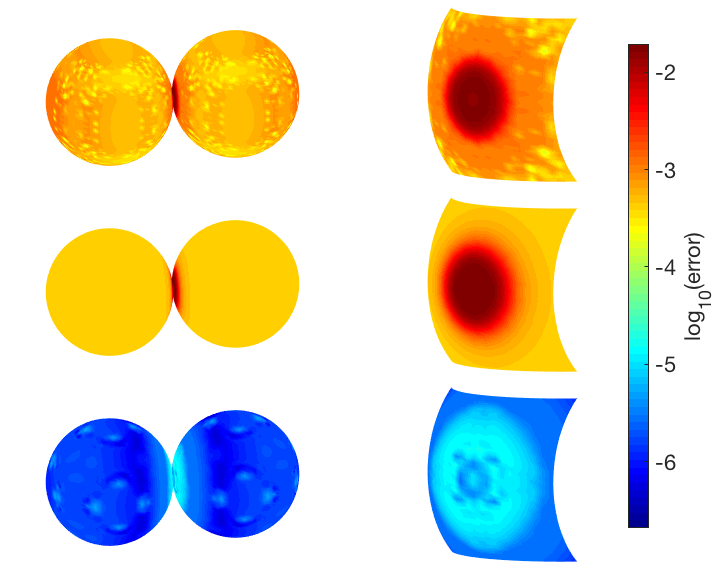}
\caption{Two spheres, solution of the integral equation \eqref{IntEq-2Surf}. Error (log scale) on the surface for $h=1/32$. Direct (top), uncorrected (middle), and corrected (bottom) solution, with the close-ups of the near singular region (right).}
\label{2Spheres}
\end{figure}

In the numerical simulation above, the integral equation is solved for a single interface, so the high order regularization of Sec. 3.3 is used and corrections are not present (unless one wishes to compute the flow off the surface). However, in the case of two or more surfaces that get close to each other, the error will deteriorate due to the near singularity issue and the high order regularization alone will not help. To demonstrate the importance of corrections, we modify the previous test in the following way. We use two unit spheres, one centered at the origin and the other at $(2,0,\epsilon)$, where $\epsilon = 1/16^3$. The integral equation is similar to \eqref{IntEq-1surf} but with a sum of two single layer and two double layer integrals, one for each interface:
\begin{align}
	\label{IntEq-2Surf}
	(\lambda_p+1) u_i(\bd{x}_0) = &-\frac{1}{4\pi \mu_0}\sum_{m=1}^{2}\int_{\pa\Omega_m} S_{ij}(\bd{x}_0,\bd{x}) [f]_j(\bd{x})dS(\bd{x}) \nonumber\\
	&+ \sum_{m=1}^{2}\frac{\lambda_m-1}{4\pi}\int_{\pa\Omega_m} T_{ijk} (\bd{x}_0,\bd{x}) u_j(\bd{x}) n_k(\bd{x})dS(\bd{x})
\end{align}
for $\bd{x}_0\in \pa\Omega_p$, $p=1,2$, $\lambda_p=\mu_p/\mu_0$. We take $\mu_0=1$ and $\lambda_1=\lambda_2=2$ as the viscosity ratio for each interface, and define the surface force as before using the surface tension $\gamma = 1+(x_1-x_c)^2$, where $x_c$ is the $x$-coordinate of the center of the sphere. We again perform successive evaluations to compute the interface velocities, similar to \eqref{Successive}. Regularization $\delta/h=2$ was used off the surface, while on the surface $\delta/h=3$ was used along with the high order regularization, as in the previous test. We solve the equation for $h=1/16, 1/32, 1/64$, and estimate the convergence rates using the error definition of \eqref{error_h}. Table~\ref{BIE-2Surf} compares the errors for three solutions: direct, uncorrected, and corrected. In the direct solution, regularization \eqref{s1}-\eqref{s3} was used without corrections. The uncorrected solution was computed using the higher regularization \eqref{s1_high}-\eqref{s3_high} for same-surface integrals, i.e., when $p=m$ in \eqref{IntEq-2Surf}, but without corrections for the other-surface case, i.e., when $p\neq m$ in \eqref{IntEq-2Surf}. Finally, the corrected solution uses higher regularization for same-surface integrals, and corrections \eqref{Corr_SL}-\eqref{Corr_DL} for the other-surface integrals. The error reported in Table~\ref{BIE-2Surf} is over the points on one sphere - it is essentially the same for the other sphere. Not only the error decreases in magnitude with corrections, but the convergence rate improves from first to third order. Figure~\ref{2Spheres} shows the error distribution (on a $\log$ scale) for the three cases. It is clear that the largest error is where the surfaces are near, and it does not improve overall when higher regularization is used for the same-surface integrals, unless the corrections are added for the other-surface integrals.

\begin{table}[htb]
\centering
\begin{tabular}{c||c|c||c|c||c|c}
\hline
\multicolumn{1}{c||}{} &
\multicolumn{2}{c||}{Direct} &
\multicolumn{2}{c||}{Uncorrected} &
\multicolumn{2}{c}{Corrected} \\
\hline
$h$ & $\lVert e_h\rVert_\infty$ & $\lVert e_h\rVert_2$ & $\lVert e_h\rVert_\infty$ & $\lVert e_h\rVert_2$ & $\lVert e_h\rVert_\infty$ & $\lVert e_h\rVert_2$ \\
\hline
1/16 & 3.53e-02 & 5.38e-03 & 3.82e-02 & 5.85e-03    & 2.34e-04 & 6.41e-05 \\
\hline
1/32 & 1.80e-02 & 1.90e-03 & 1.96e-02 & 2.03e-03    & 2.66e-05 & 3.36e-06 \\
\hline
\end{tabular}
\caption{Two spheres, error in the solution of the integral equation \eqref{IntEq-2Surf}. Direct: regularization \eqref{s1}-\eqref{s3} used everywhere without corrections. Uncorrected: still no corrections but higher regularization \eqref{s1_high}-\eqref{s3_high} used in same-surface integration ($p=m$). Corrected: higher regularization in same-surface integration, and corrections \eqref{Corr_SL}-\eqref{Corr_DL} used for other-surface integration ($p\neq m$). Regularization $\del/h=3$ for same-surface and $\del/h=2$ for other-surface.}
\label{BIE-2Surf}
\end{table}

\section{Conclusions}

The numerical results we have performed in this paper are in agreement with the analytical prediction of uniformly third order spatial convergence for the computation of single and double layer integrals of the form \eqref{SingleLayer} and \eqref{DoubleLayer}. This is true when the evaluation point is near the surface, as is the case when two interfaces are close to each other. The accurate solution is obtained by regularizing the kernels and adding analytically derived correction terms to eliminate the first and second order regularization error terms. When the evaluation point is on the surface, a much higher accuracy is achieved without corrections, by improving the way the integrands are regularized. For this case, it might help to derive correction terms for the error due to the discretization of the integrals, although with an appropriate choice of the regularization parameter, this seems rarely necessary.

The error decays rapidly away from the surface, but our numerical results suggest that it might still be somewhat larger when the evaluation point is very close to the surface. One might experiment then with an interpolation technique such as \cite{ying06}, where the value very near the surface is interpolated from corrected values further away, rather than computed directly with corrections.

One of the advantages of the corrections method is that it does not increase the computational complexity of the overall method. Specifically, with $N$ quadrature points and $M$ evaluation points, computing the integrals will require $O(NM)$ CPU time, while the corrections add $O(M)$ to this. The computational efficiency of the algorithm can be improved then irrespective of the corrections, by using a fast summation algorithm such as a treecode \cite{treecode} or a fast multipole method \cite{tornberg08}. Such methods decrease the CPU time generally to $O(M\log N)$ or $O(M)$, respectively. 


\bibliographystyle{plain}

\end{document}